\pdfoutput=1
\documentclass[fontsize=11pt,a4paper,DIV=12]{scrartcl}

\usepackage[T1]{fontenc}

\usepackage{enumerate,amsthm, amsmath, amsfonts, amssymb, mathrsfs, complexity, wasysym, graphicx,colortbl, url, hyperref, hypcap, shuffle, xargs, multicol, overpic, pdflscape, multirow, hvfloat, minibox, accents, array, multido, xifthen, xspace, ae, aecompl, blkarray, pifont, mathtools, etoolbox, dsfont, verbatim, stackrel, stmaryrd, xcolor, libertine, dynkin-diagrams, thmtools, tikz-cd, xfrac, nicefrac}

\usepackage{shapepar, microtype, lipsum}
\usepackage{pgfplots}
\pgfplotsset{compat=1.9}
\usepackage{tikz}
\usetikzlibrary{calc,fit,intersections,folding}
\usepackage{pstricks-add}
\usetikzlibrary{arrows.meta,angles,arrows,quotes,backgrounds}

\hypersetup{colorlinks=true, citecolor=darkblue, linkcolor=darkblue,urlcolor=darkblue}
\usepackage[noabbrev,capitalise]{cleveref}

\usepackage[subrefformat=simple,labelformat=simple]{subcaption}

\usepackage[backend=biber,style=alphabetic,sorting=ynt, maxbibnames=99, maxcitenames=4]{biblatex}
\hypersetup{colorlinks=true, citecolor=darkblue, linkcolor=darkblue,urlcolor=darkblue}

\usepackage[subrefformat=simple,labelformat=simple]{subcaption}

\definecolor{darkblue}{rgb}{0,0,0.7} % darkblue color
\definecolor{darkred}{rgb}{0.7,0,0} % darkred color
\definecolor{green}{RGB}{57,181,74} % darkgreen color
\definecolor{violet}{RGB}{147,39,143} % violet color

\definecolor{color1}{HTML}{6500ca}
\definecolor{color2}{HTML}{447bcc}
\definecolor{color3}{HTML}{82b96f}
\definecolor{color4}{HTML}{dcaa3c}
\definecolor{color5}{HTML}{db2021}
\definecolor{color6}{HTML}{d64998}

\newcommand{\darkblue}{\color{darkblue}} % darkblue
\newcommand{\defn}[1]{\textsl{\darkblue #1}} % emphasis of a definition

\DeclareMathOperator{\N}{\mathbb{N}}
\let\A\relax
\DeclareMathOperator{\A}{\mathcal{A}}

\newcommand{\setclr}[1]
{
    \if#11\def\clr{color1}\fi
    \if#12\def\clr{color2}\fi
    \if#13\def\clr{color3}\fi
    \if#14\def\clr{color4}\fi                   
    \if#15\def\clr{color5}\fi
}

\newcommand{\cupdot}{\mathbin{\mathaccent\cdot\cup}}
\newcommand\bigcupdot{\mathop{\mathaccent\cdot{\bigcup}}}

\newtheorem{theorem}{Theorem}
\newtheorem{observation}{Observation}

\newtheorem{remark}{Remark}

\newtheorem{lemma}{Lemma}
\newtheorem{proposition}{Proposition}
\newtheorem{corollary}{Corollary}
\newtheorem*{theorem*}{Theorem}

\newtheorem{example}{Example}

\DeclareTOCStyleEntry[
    indent=0em,
    pagenumberformat=\normalfont,
    beforeskip=1pt plus .2pt,
    entryformat=\normalfont
    ]{tocline}{section} 
    
\begin{document}

\title{\huge Balanced Gray Codes for Permutations and
	Rainbow Cycles for Associahedra}
\author{
	Robert Lauff\,\thanks{Technische Universit\"at Berlin, Germany, lauff@math.tu-berlin.de} \and
	Lucca Tiemens\,\thanks{No affiliation, lucca@tiemens.de}
}
\date{}

\maketitle
\thispagestyle{empty}

\begin{abstract}
	We settle the problem of constructing a balanced transposition Gray code for permutations of $[n] := \{1, \dots, n\}$ with $n \in \mathbb{N}\setminus\{0\}$. More generally, we obtain a~$2(m-2)!$-rainbow cycle for the permutations of $[n]$ for $m \in [n]$, a notion recently introduced by \citeauthor[]{Felsner2020}. Furthermore, we extend a result of theirs by presenting a $r$-rainbow cycle for the classical associahedron $\mathcal{A}_{n}$ for $r \sim 2^n$. Previously, the longest rainbow cycle for the associahedron was linear in length.

	For even $n$, we also construct a balanced Gray code for permutations of $[n]$, using only cyclically adjacent transpositions, complementing the construction for odd $n$ by \citeauthor[]{Gregor2023}.
	Additionally, we show that the Permutahedron $P_{n}$ admits a $2$-rainbow cycle for all $n\ge5$ and a $3$-rainbow cycle for odd $n\ge3$.

\end{abstract}

\pagebreak
\thispagestyle{empty}
\setcounter{tocdepth}{2}
\tableofcontents

\section{Introduction}
\label{sec:intro}

\setcounter{page}{1}

In the study of combinatorial generation, one is given a class of combinatorial objects and the task of generating them exhaustively. A \emph{combinatorial Gray code} is such an enumeration with the additional properties that two neighboring objects differ only by a local change and no object is listed twice. Many known Gray codes lead to elegant and efficient algorithms for solving the related combinatorial generation task. Because only local changes are performed in each step, these algorithms can be very memory efficient. Gray codes are named after the physicist Frank Gray who worked at Bell Labs, where he patented a method, today known as \emph{binary reflected Gray code} \cite{Knuth2011}, to generate all $2^{n}$ bit strings of length $n$ \cite{Gray1953}. In his construction two consecutive elements differ only by flipping a single bit. Since then, Gray codes have been found for several combinatorial structures including permutations, set partitions, integer partitions, spanning trees, Hamilton cycles, etc. (see \cite{Muetze2023} or \cite{Sav97} for surveys on the subject). In practice, a Gray code is a Hamilton path (or cycle) in a \defn{flip graph} of the combinatorial structure, defined by having the objects as vertices and the local changes as edges. For example, the binary reflected Gray code for length $n$ bit strings is a Hamilton cycle in the~$n$-dimensional cube and the local changes are bit flips.

When constructing Gray codes, there are desirable additional properties to aim for, for example controlling the transition count for the different kinds of flips available. A Gray code with a perfect balance of flip operations is called a \defn{balanced Gray code}. The study of Gray codes for bit strings revealed that it is possible to design such a code where the total flip counts in each of the $n$ coordinates differ by at most $2$ \cite{bakos1968, Bhat1996}. When $n$ is a power of $2$ it is possible that each bit is flipped exactly~$\frac{2^{n}}{n}$ times \cite{wagner1991}. This raises the question if such a balancing of changes is also possible for other classes of combinatorial objects. Specifically for permutations, Mütze has asked about the existence of a balanced Gray code using transpositions as flips in \cite{Muetze2023} motivated by \cite{Felsner2020}. We are able to give a positive answer. On the other hand, as initial progress towards a balanced Gray code for the permutahedron, the flip graph of permutations having only adjacent transpositions as edges, \citeauthor{Gregor2023} created a balanced Gray code using the cyclically adjacent transpositions for odd $n$ \cite{Gregor2023}. We now also provide the construction for even $n$.

It is worth noting that there are many known non-balanced Gray codes for permutations (see for example \cite{Lipski1979} , \cite{Kompel’makher1976} and \cite{slater1978} or \cite[chapter 5]{Muetze2023} for a detailed overview). A famous procedure to generate all elements of the symmetric group by applying only adjacent transpositions is called the \emph{Steinhaus-Johnson-Trotter (SJT)} algorithm \cite{Johnson1963}, \cite{Trotter1962} or method of \emph{plain changes} \cite[Sec. 7.2.1.2]{Knuth2011}. An interesting advancement of the SJT algorithm is the \emph{doubly-adjacent} Gray code by \citeauthor{Compton1993} given in \cite{Compton1993} where any two consecutive adjacent transpositions~$(i, i+1)$ and $(j, j+1)$ additionally satisfy $|i-j| = 1$.

Balanced gray codes are generalized by \defn{$r$-rainbow cycles}, which are cycles in flip graphs using each kind of flip, also called color, exactly $r$ times. The color of the flip has a combinatorial interpretation in practice. Rainbow cycles are known to exist in a variety of flip graphs. \citeauthor[]{Felsner2020} \cite{Felsner2020} constructed 1- and 2- rainbow cycles on the classical associahedron, $1,\ldots,(|X|-2)$-rainbow cycles on the flip graph of plane spanning trees on point set $X$ in general position, 1- and 2-rainbow cycles in the flip graph of non-crossing perfect matchings on $2m$ points in convex position (1-rainbow for $m\in\{2,4\}$, 2-rainbow for $m\in\{6,8\}$) and 1-rainbow cycles for $k$-subsets of $[n]$ under element exchange for $2\le k<\frac{n}{3}$, as well as some negative results. For $2$-subsets of $[n]$, if $n$ is odd, they find two edge disjoint 1-rainbow Hamilton cycles, i.e. two balanced Gray codes. Last but not least, they construct 1-rainbow cycles for permutations of $[n]$ with transpositions as flips, for even $\lfloor\frac{n}{2}\rfloor$. Our results extend their work on the classical associahedron as well as on permutations. Rainbow cycles on the associahedron are especially interesting, because they relate to so called \emph{facet-Hamiltonian cycles} \cite{facethamiltonicity} (SODA25). Details about this connection are provided in \Cref{sec:results and preliminaries}.

Outside the flip graph setting, rainbow cycles and paths have been studied much earlier. A well-known conjecture by \citeauthor[]{Andersen89} \cite{Andersen89} states that the complete graph on $n$ vertices, with any proper edge coloring, contains a path of length $n-2$ with distinct colors along its edges. See \cite{Alon16} by \citeauthor[]{Alon16} and \cite{Balogh17} by \citeauthor[]{Balogh17} for progress towards this conjecture.

\subsection*{Outline of the paper}

The paper is structured as follows. In \Cref{sec:results and preliminaries} we present the main results in more detail, giving basic definitions needed for their formulation and present proof sketches. In \Cref{sec:associahedron} we present the construction of $r$-rainbow cycles for $r\sim 2^n$ in the associahedron, where~$n$ is the number of vertices of the triangulated polygon, see \Cref{thm:associahedron}. In \Cref{sec:balancedOdd} and \ref{sec:balancedEven} we consider the flip graph of permutations of $[n]$ where all transpositions are admissible and obtain $2(m-2)!$-rainbow cycles for $m\in[n]$, see \Cref{thm:perm_rainbow}. In \Cref{sec:cyclicallyBalanced} we restrict the admissible transpositions to cyclically adjacent ones and construct a balanced Gray code for even $n$, see \Cref{thm:cyclicallyballanced}. \Cref{sec:permutahedron} covers the construction for 2- and 3-rainbow cycles in the permutahedron, see \Cref{thm: 2 and 3 rainbow for perm}. Finally, in \Cref{sec:conclusion}, we discuss the results and consider open problems.

\section{Results and preliminaries}\label{sec:results and preliminaries}

In this paper we settle the problem of constructing a balanced transposition Gray code for permutations, with all transpositions and restricted to cyclically adjacent transpositions, and we construct $r$-rainbow cycles on the classical associahedron for $r\sim 2^n$, where $n$ is the size of the triangulated polygon. More specifically, in the flip graph of permutations of $[n]$ with all transpositions allowed, we construct $2(m-2)!$-rainbow cycles for $m\in[n]$. The maximum value $m=n$ gives the balanced Gray code. Furthermore, in the case that only cyclically adjacent transpositions are considered, we construct a balanced Gray code for even $n$. This complements a construction by Gregor, Merino and Mütze \cite{Gregor2023} who proved the result for odd $n$. Restricting further to just adjacent transpositions, we present a $2$ and $3$-rainbow cycle for the permutahedron.

The question about the existence of a balanced Gray code using all transpositions appears as $P30$ in the survey \citetitle[p. 39]{Muetze2023} by \citeauthor[]{Muetze2023} \cite{Muetze2023}. The original source for this problem is the before-mentioned paper \citetitle[]{Felsner2020} by \citeauthor[]{Felsner2020} \cite{Felsner2020}, which asks, in particular, for which $r$ it is possible to create $r$-rainbow cycles for permutations. Hence, our result also advances the original question. 
Furthermore, they construct $1$ and $2$-rainbow cycles in the classical associahedron, which we are able to extend to exponential $r$.

In the following parts, we first discuss the preliminaries and give proof sketches.

\subsection{The Associahedron}

\subsubsection{Preliminaries}

For a given $n\in\N$, the vertices of the \defn{classical associahedron} are the triangulations of a convex $n$-gon. Let $T$ be some triangulation and let some 4-gon $g$ in it be given. Then exactly one of the two diagonals of~$g$ is in $T$. The \defn{flip} of this diagonal replaces it by the other diagonal of $g$. The classical associahedron has as edges the pairs of triangulations which differ by a flip and is denoted by $\A_n$. We color the edges by the diagonal which is flipped. Hence we seek cycles in $\A_n$ which flip each diagonal a fixed number $r$ times. A 1-rainbow cycle in $\A_n$ is also known as a facet-Hamiltonian cycle. A facet-Hamiltonian cycle on a polytope $P$ is a cycle $C$ on the skeleton of $P$ such that $C\cap f$ is non-empty and connected for any facet $f$ of $P$. The facets of the associahedron are given by those sets of triangulations sharing a given diagonal. Hence, if each diagonal is flipped exactly once along a cycle $C$, then the triangulations containing it form a connected arc of $C$. A $r$-rainbow cycle is hence a cycle $C$ on $\A_n$ which visits each of its facets $f$ exactly $r$ times, i.e. $C\cap f$ has $r$ connected components. We also define \defn{r-rainbow paths} to be paths on $\A_n$ which flip every diagonal of the $n$-gon exactly $r$ times, but where we also count the diagonals of the final triangulation.

\begin{figure}[hbt]
        \centering
            \centering
            \scalebox{0.3}{
\begin{tikzpicture}[rotate = 12.86]
    %blue points
    \foreach \i in {1,...,15} {
        \node[fill, circle] (\i) at (360*\i/15:5) {};
    }
    \begin{scope}[on background layer]
        \foreach \i in {1,...,14} {
            \pgfmathparse{int(\i+1)}
            \draw[gray] (\i.center) -- (\pgfmathresult.center);
        }
        \draw[gray] (15.center) -- (1.center);
    \end{scope}

    \foreach \x/\y in 
    {1/2, 1/3, 1/4, 1/5, 1/6, 6/15, 6/14, 6/13, 13/7, 7/12, 7/11, 7/10, 10/8}
    {
        \draw[thick] (\x.center) -- (\y.center);
    }
\end{tikzpicture}
}
        \caption{A zigzag triangulation}
        \label{fig:zigzag}
    \end{figure}

A \defn{zigzag-triangulation} is a triangulation of the $n$-gon whose weak dual is a path, see \Cref{fig:zigzag}. This path will be referred to as the \defn{dual path}. Equivalently, every triangle in a zigzag-triangulation has an edge incident to the outer face. Such a triangulation has exactly two vertices of degree $2$, called the \defn{endpoints}. A \defn{star} is a triangulation of the $n$-gon such that there is a vertex of degree $n-1$. This vertex is called the \defn{center}.
Given a triangulation~$T$ in $\A_n$ and a permutation $\pi$ of its diagonals, the triangulation obtained by flipping the diagonals of $T$ in the order of~$\pi$ is denoted as \defn{$T(\pi)$}. The resulting path is denoted as \defn{$T\xrightarrow{\pi}T(\pi)$}. Note that the new diagonals introduced along this path are exactly the diagonals present in $T(\pi)$. We say that $T_1$ and $T_2$ are \defn{permutation-adjacent} if there is a permutation $\pi$ such that $T_2=T_1(\pi)$. The diagonals of $T(\pi)$ are in bijection to the diagonals of $T$, where each diagonal is associated to its origin in $T$. Hence, $\pi$ induces a natural permutation of the diagonals of $T(\pi)$. Flipping the diagonals of $T(\pi)$ in that order yields the triangulation~\defn{$T(\pi^2)$}, and similarly we define \defn{$T(\pi^k)$}, $k\in\N$. The path obtained is denoted by \defn{$T\xrightarrow{\pi^k}T(\pi^k)$}.

The point of permutation-adjacency is the following lemma, which follows directly from the definition:
\begin{lemma}
    Let $C$ be a cycle on $\A_n$ which splits into parts as follows:
    \[
        C=T_1\xrightarrow{\pi_1}T_2\xrightarrow{\pi_2}\ldots\xrightarrow{\pi_{k-1}}T_k\xrightarrow{\pi_k}T_1
    \]
    for triangulations $T_1,\ldots,T_k$ and permutations $\pi_1,\ldots,\pi_k$. Then a diagonal of the $n$-gon is flipped along a traversal of $C$ exactly as often as it appears in the triangulations $T_1,\ldots,T_k$.
\end{lemma}

Our goal will hence be to find a suitable set of triangulations which are permutation adjacent so that they form a cycle, and such that each diagonal of the $n$-gon appears $r$ times among them. These will be subsets of \defn{imbalanced zigzag triangulations}. Loosely speaking, they consist of any zigzag triangulation, called the label part, and a star glued together. A proper definition is given in \Cref{sec:associahedron}.

\begin{figure}[hbt]
    \centering
    \include{Pictures/imbalanced_example}
    \caption{Imbalanced zigzag triangulations. On the left the case where $n$ is even. In this case the number of diagonals is odd. The longest diagonal then counts neither to the star nor the label part. On the right the case where $n$ is odd. In both cases the star is blue and the label part is red.}
    \label{fig:imbalanced}
\end{figure}

\begin{remark}
    The classical associahedron has many more combinatorial interpretations. For example, its vertices can be viewed as the set of rooted binary trees with a fixed number of vertices (where each non-root vertex is either the left or right child of its parent, even if the parent has only one child). Then the edges correspond to tree rotations. The bijection is obtained by considering the dual graph of a triangulation, but remembering the side of the triangle which needs to be passed to get to the next vertex of the dual graph. Zigzag triangulations are thus in bijection to rooted paths.
\end{remark}

\begin{restatable}{theorem}{associahedron}\label{thm:associahedron}
    Let $1\le k\le\frac{2^{m-1}}{\binom{m-1}{6}(2^6-1)}-2$, where $m=\lceil\frac{n-4}{2}\rceil$. There is a $4k$-rainbow cycle on $\A_n$.
\end{restatable}
\begin{proof}[Proof sketch]
    See \Cref{fig:construction imbalanced} for an illustration of the construction.
    
    The 2-rainbow cycles for the associahedron from \cite{Felsner2020} are constructed by starting with a star. This star is rotated around by flipping the edges in the cyclic order around its center. It is shown that the result of flipping each edge in this order is another star, which is rotated by one vertex. Then, proving that each diagonal of the $n$-gon appears in exactly two stars finishes the proof.
    
    We replace the stars with imbalanced zigzag triangulations and show that such a triangulation can also be rotated by one vertex by flipping all edges in a particular order (\Cref{lem:rotation}). In fact, any zigzag triangulation is permutation adjacent to its rotation by one vertex (in either direction). After showing that each edge of the $n$-gon appears in exactly two triangulations from the family of rotations of a zigzag triangulation (\Cref{lem:every edge twice}), we have a way of constructing many different 2-rainbow cycles (\Cref{lem:2rainbow}). Then we can glue these together in the following way: We show that any zigzag triangulation is permutation adjacent to a star (\Cref{lem:adj to stars}). This allows us to transform one imbalanced zigzag triangulation into another. We first transform the star of the one into the label of the other, and then the label of the one into the star of the other.

    The stars allow for easy identification and assure that no triangulation is repeated along the so defined cycle.
\end{proof}

\newcommand{\Pall}{\mathcal{P}^{\textit{all}}}
\newcommand{\Pcadj}{\mathcal{P}^{\textit{cadj}}}
\newcommand{\Padj}{\mathcal{P}^{\textit{adj}}}

\subsection{Permutation Gray Codes}

We now come to our results about permutations. We start with the needed preliminaries. Let \defn{$S_{n}$} denote the set of permutations of $[n] := \{1,\dotsc,n\}$. Throughout the paper write ${\sigma \in S_{n}}$ as a linear array $\sigma_{1}\dotsc \sigma_{n}$ defined by $\sigma(i) = \sigma_{i}$. Further, we denote a \defn{transposition} by $(i, j)$ and write $\sigma \xrightarrow{(i, j)} \sigma'$ to mean that $\sigma'$ differs from $\sigma$ by interchanging the numbers (not the indices) $i$ and $j$ in the linear array representation. For example, ${213 \xrightarrow{(1, 3)} 231}$. A transposition $(i, j)$ is called \defn{adjacent} if $|j - i| = 1$. Otherwise, the transposition is called \defn{wide}. The \defn{cyclically adjacent} transpositions on $S_{n}$ are all adjacent transpositions together with $(1, n)$.

A \defn{Gray code for $S_{n}$} is a listing $(\pi_{i})_{i\in[n!]}$ of all elements of $S_{n}$, such that ${\pi_{i} \xrightarrow{t_{i}} \pi_{i+1}}$ for some transposition~$t_{i}$. Equivalently, a Gray code can be represented by listing just the transpositions $(t_i)_{i \in [n!-1]}$. We will represent permutation Gray codes in this form throughout the paper.
To this end, we sometimes write \defn{$\pi \xrightarrow{G} \pi'$} to denote a Gray code $G$ starting at $\pi$ and ending at $\pi'$. Note that if no starting permutation is given, one can start a Gray code represented by transpositions on any permutation, as long as all transpositions are well-defined (meaning all numbers appear in the permutation). For a Gray code $G$ we denote by \defn{$\overleftarrow{G}$} the reverse Gray code, defined by reversing the list of transpositions.

In the following we denote by $\Pall_n, \Pcadj_n$ and $\Padj_n$ the following flip graphs with $S_n$ as vertex set. In each of these, we connect permutations who differ by an admissible transposition. In $\Pall_n$, all transpositions are admissible, in $\Pcadj_n$ the cyclically adjacent transpositions are admissible and in $\Padj_n$ only the adjacent transpositions are admissible. The edges are colored by the transposition applied.

Finally, we present the definition that enables our construction on $\Pall_n$.
An \defn{almost balanced Gray code} $G$ for $\Pall_{n}$ is a Gray code such that all wide transpositions appear exactly $2(n-2)!$ times in $G$.
Equivalently, a Gray code $G$ for $\Pall_n$ is almost balanced if only adjacent transpositions can appear less or more often then $2(n-2)!$ times in $G$.

We obtain the following results.

\begin{restatable}{theorem}{permrainbow}\label{thm:perm_rainbow}
    For $n\in \mathbb{N}\setminus\{0\}$ and $m \in [n]$, $\Pall_n$ admits a $2(m-2)!$-rainbow cycle. In particular, a balanced transposition Gray code is obtained by setting $m = n$.
\end{restatable}
\begin{proof}[Proof sketch]
    The proof is done in $3$ stages. First, we construct a balanced Gray code for $\Pall_n$ for odd~$n\in \mathbb{N}$. Then the balanced Gray code for even $n\in \mathbb{N}$ is presented. Finally, we show that a balanced Gray code for $\Pall_n$ for any $n$ can be lifted to a $2(n-2)!$-rainbow cycle in $\Pall_m$ for all $m>n$, which completes the proof.
    We now sketch each stage.

    Let $n\in \mathbb{N}$ be odd. We obtain the balanced Gray code for $\Pall_n$ inductively. The code \[G_3 := ((1,3), (1,2), (2,3), (1,3), (1,2))\] for $\Pall_3$ is our base case. We have $123 \xrightarrow{G_3} 132$ such that $(2,3)$ closes the cycle.
    Suppose we are given a balanced Gray code $G_n$ for odd $\Pall_n$ such that $(n-1, n)$ closes the cycle. We then transform $G_n$ into another balanced Gray code $H_n$ for $\Pall_n$ such that now $(1,2)$ closes the cycle (see \Cref{lem:balancedCodeH}). We proceed by showing that for the permutation $\omega_k := (k \dotsc n+1)$ the code $L$ defined by
    \[
        H_n, (n,n+1), \overleftarrow{\omega_n(H_n)}, (n-1, n), \omega_{n-1}(H_n), \dotsc , \overleftarrow{\omega_3(H_n)}, (2, 3), \omega_2(G_n), (1, 2), \overleftarrow{\omega_1(G_n)}
    \]
    is an almost balanced Gray code for $\Pall_{n+1}$ such that \[123\dotsc n(n+1) \xrightarrow{L} 234\dotsc n(n+1)1\] (see \Cref{prop:almostbalancedCode}). This left-rotation encoded by $L$ is then utilized to prove that for the cyclic rotation $\sigma := (1\dotsc n+2)$ the code $G_{n+2}$ given by
    \[
        L, (1,n+2), \sigma(L), (1,2), \sigma^2(L), \dotsc , \sigma^n(L)
    \]
    is again a balanced Gray code for $\Pall_{n+2}$ such that $(n+1, n+2)$ closes the cycle (see \Cref{prop:balancedCode}). The intuition is that the repeated rotations cause the unbalanced adjacent transpositions of $L$ to balance as the cyclically adjacent transpositions in $G_{n+2}$, while the already balanced wide transpositions of $L$ simply lift into balanced appearances in $G_{n+2}$, since $L$ is by construction almost balanced.

    In the next stage let $n>0$ be even. Similar to the construction of $L$ in the previous stage, we additionally construct an almost balanced Gray code $R$ for $\Pall_n$ such that \[123\dotsc n \xrightarrow{R} n123\dotsc n-1\] yielding a right-rotation with the additional property that the multiset of transpositions of $L$ and $R$ is the same (see \Cref{prop:almostbalancedCodeRight}). We then create a non-cyclic Gray code $H$ for $\Pall_{n+1}$ using $\sigma := (1 \dotsc n+1)$ by
    \[
        \overleftarrow{L}, (1,n+1), \overleftarrow{\sigma(R)}, (1,2), \overleftarrow{\sigma^2(R)}, \dotsc , \overleftarrow{\sigma^{n-1}(R)}
    \]
    which has the same multiset of transpositions as the balanced Gray code for $\Pall_{n+1}$ obtained in the previous stage (\Cref{prop:secondGrayCode}), which we now denote by $G$. Let $\sigma := (n+2 \dotsc 1)$ be the reversed cyclic rotation. Then the alternation between increasingly rotated versions of $H$ and $G$ given by
    \[
        G, (n+1, n+2), \sigma(H), (n,n+1), \sigma^2(G), (n-1, n), \sigma^3(H),\dotsc,\sigma^n(G),(1,2), \sigma^{n+1}(H)
    \]
    is a balanced Gray code for $\Pall_{n+2}$ (see \Cref{prop:balancedEven}). The intuition is that, since $G$ and $H$ have the same multiset of transpositions and $G$ is balanced for $\Pall_{n+1}$, the repeated rotations $\sigma$ cause each transposition in $G$ to lift to balanced appearances in $\Pall_{n+2}$.

    The proof is then finished by observing that a rainbow cycle for $\Pall_n$, $n\in \mathbb{N}$, can be lifted to $\Pall_{n+1}$ by pushing $n+1$ from the left to the right and back on an appropriate number of permutations while traversing the cycle (see \hyperref[proof:theorem2]{proof} of \Cref{thm:perm_rainbow}).
\end{proof}

Furthermore, our method for constructing this balanced Gray code enabled us to complement a result by \citeauthor[]{Gregor2023}. In \cite{Gregor2023} the authors present for odd~$n \in \mathbb{N}$ a balanced Gray code for $\Pcadj_n$. We now obtained this result for even $n$ as well.

\input{Examples/cadj-and-all-4}

\begin{restatable}{theorem}{cyclicallyBalanced}\label{thm:cyclicallyballanced}
  Let $n \in \mathbb{N} \setminus {0}$ be even. Then there is a balanced Gray code for $\Pcadj_n$, meaning that each of the $n$ cyclically adjacent transpositions appears exactly $\frac{n!}{n} = (n-1)!$ times within the cycle.
\end{restatable}
\begin{proof}[Proof sketch]

    The idea is similar to the construction of the balanced Gray code for $\Pall_n$ in the even case. Contrary to \Cref{thm:perm_rainbow}, the transpositions in the proof of \Cref{thm:cyclicallyballanced} are applied to the indices and not the entries of a permutation.

    Let $n\in \mathbb{N}\setminus \{0\}$ be even. By \cite{Tchuente1982} there is a Hamilton path between any permutation $\pi_1\pi_2\dotsc\pi_n$ and $\pi_2\pi_3\dotsc\pi_n\pi_1$ on the Permutahedron. We denote the sequence of transpositions from this path by $L$, while we denote the reversed sequence by $R$. Using $L$ and $R$ we construct two Gray codes for $\Padj_{n+1}$. Let $\sigma := (12\dotsc n+1)$. The first one, call it $G$, is given by 
    \[
        \sigma(L), (1,2), \sigma(L), (1,2),\dotsc,\sigma(L)
    \]
    and the second one, denoted by $H$, is given by
    \[
        \sigma(R), (1,2), \sigma(L), (1,2),\dotsc,\sigma(L).
    \]
    For those codes it holds that 
    \[ 
        12\dotsc(n+1) \xrightarrow{G} 213\dotsc(n+1) \text{ (see \Cref{prop:cyclical_G}) }
    \] 
    and 
    \[ 
        13\dotsc(n+1)2 \xrightarrow{H} (n+1)12\dotsc n \text{ (see \Cref{prop:cyclical_H}). }
    \]
    In particular, since $L$ and $R$ have the same multiset of transpositions, $G$ and $H$ have the same multiset of transpositions. We then prove that the Gray code obtained by alternating between increasingly rotated versions of $G$ and $H$ given by 
    \[
        G, (1,n+2),\sigma(H),(1,2),\sigma^2(G),(2,3),\sigma^3(H),\dotsc,\sigma^{n+1}(H)
    \]
    is balanced for $\Pcadj_{n+2}$, where $\sigma := (12 \dotsc n+2)$ (see \hyperref[proof:theorem3]{proof} of \Cref{thm:cyclicallyballanced}). The intuition is that each adjacent transposition from $G$ gives rise to all cyclically adjacent transpositions through the rotations. Since the transpositions in $G$ and $H$ are the same, this leads to fully balanced appearances in the resulting code.

\end{proof}

Lastly, we investigate a classic flip graph that was not yet considered in this context, the permutahedron $\Padj_n$. In \cite{facethamiltonicity} a so called facet-Hamiltonian cycle is found for the permutahedron. This is equivalent to a rainbow-cycle of the permutahedron where the flips are colored by the new prefix they introduce. For example the flip $(3,4)$ introduces the new prefix $\{2,5,6\}$ in $5,2,3,6,4,1\to 5,2,6,3,4,1$. However, here we color by the adjacent transposition applied.

\begin{restatable}{theorem}{permutahedron}\label{thm: 2 and 3 rainbow for perm}
    The permutahedron $\Padj_n$ admits a 2-rainbow cycle for all $n\ge 5$ and a 3-rainbow cycle for all odd $n\ge 3$.
\end{restatable}
\begin{proof}[Proof sketch]
    The proof of this theorem is a simple induction on $n$.
\end{proof}

\input{Examples/Pall5}

\pagebreak

\section{Rainbow Cycles on the Associahedron}\label{sec:associahedron}

Let us fill in all the details left out in the proof sketch.

\begin{lemma}\label{lem:adj to stars}
    Let $T$ be a zigzag triangulation, $i$ a point of the $n$-gon and $r$ a ray shot from $i$ in some direction. Flipping the diagonals of $T$ hit by $r$ in the order that they are hit, causes them all to be incident to $i$.
\end{lemma}
\begin{proof}
    Assume that $r$ hits some diagonal $d$ and that the segment between $i$ and $d$ along $r$ is not intersected by another diagonal. Then $d$ splits a unique 4-gon in half and this 4-gon has $i$ as one of its vertices. The other diagonal of this 4-gon is therefore incident to $i$.
\end{proof}
By choosing the ray from one endpoint to the other, we obtain the following:
\begin{corollary}\label{cor:adj to stars}
    Every zigzag triangulation $T$ is permutation-adjacent to each of the stars with center at an endpoint of $T$.
\end{corollary}

Given a triangulation $T$, let $[T]$ be the set of triangulations in $\A_n$ obtained by rotating $T$. If $T$ has no rotational symmetries, then $|[T]|=n$. We call $[T]$ the \defn{rotation family} of $T$.

\begin{lemma}\label{lem:every edge twice}
    Let $T$ be a zigzag-triangulation. Then every possible diagonal of the $n$-gon appears among the triangulations in $[T]$ exactly twice.
\end{lemma}
\begin{proof}
    Any inner edge of the $n$-gon splits the remaining vertices into two sets. Associate to an edge the pair of numbers that are the cardinalities of these sets. Note that $T$ has two edges associated to any possible pair of sizes, except $\{\frac{n-2}{2},\frac{n-2}{2}\}$ if $n$ is even. This is because the inner edge incident to the face of an endpoint is associated to $\{1,n-3\}$, the next edge along the dual path is associated to~$\{2,n-4\}$ and so on until we arrive at $\{n-3,1\}=\{1,n-3\}$ on the other side. If $n$ is even, then $T$ has exactly one edge splitting the $n$-gon in half which is associated to $\{\frac{n-2}{2},\frac{n-2}{2}\}$. Note that this edge is the only one having a rotational symmetry. Now, fix any interior edge $e$ of the $n$-gon. If it is not associated to $\{\frac{n-2}{2},\frac{n-2}{2}\}$ then there are two edges of $T$ which are associated to the same pair of numbers. Each will line up with $e$ in exactly one rotation of $T$ in $[T]$. If $e$ is associated to $\{\frac{n-2}{2},\frac{n-2}{2}\}$, then the edge of $T$ associated to $\{\frac{n-2}{2},\frac{n-2}{2}\}$ will line up with it twice along the rotation because of its rotational symmetry.
\end{proof}

\begin{figure}
        \centering
        \includegraphics[width=0.8\linewidth]{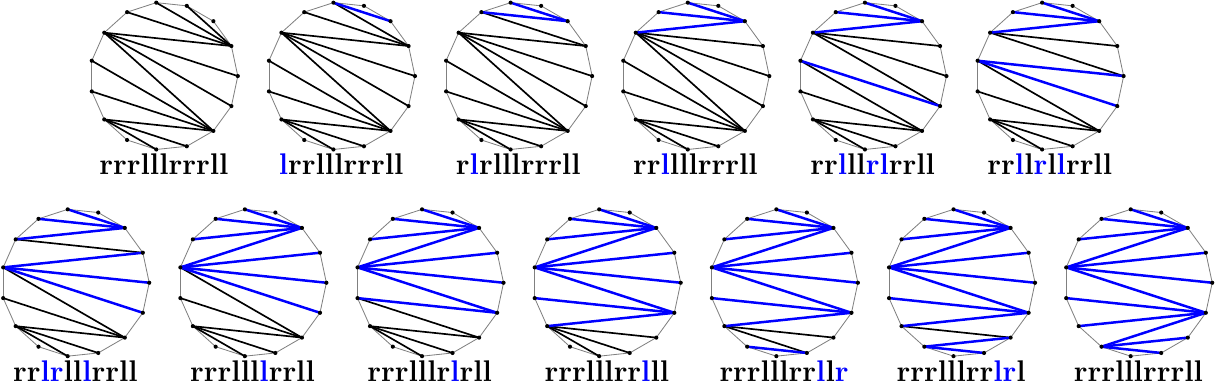}
        \caption{Rotating a zigzag triangulation. The word corresponding to the upper right endpoint is written underneath each triangulation. The positions at which the word does not match the word of the initial triangulation are marked blue.}
        \label{fig:rot and labels}
    \end{figure}

To any zigzag triangulation we associate 2 \defn{words}. Start at either endpoint and look at the edges in that order. The next edge along always includes a new vertex to the right or to the left from the perspective of walking the dual path. We save this information as a word on $\{l,r\}$, see \Cref{fig:rot and labels}. If~$i$ is the endpoint of the triangulation $T$ we started from, then $\omega_i(T)$ denotes the resulting word. Note that the words are invariant under rotation, hence we can assign the pair of words to the rotation families of the zigzag triangulations. The length of the words is $n-4$, one less than the number of diagonals. Given any word $w$ on $\{l,r\}$, we define \defn{$\bar{w}$} to be the word obtained from $w$ by reversing and exchanging~$l$'s and $r$'s. If $w$ is a word of a zigzag triangulation associated to one end point, then~$\bar{w}$ is the word associated to the other one. Given a word $w$ on $\{l,r\}$, we say that its \defn{blocks} are the longest consecutive sequences of $l$'s and $r$'s in it, for example $w=lllrrrlllrl$ consists of the blocks~$lll|rrr|lll|r|l$. We write $w=b_1b_2\ldots b_k$ where the $b_i$ are the blocks.

Next, the following lemma will give us a way to explore all elements of $[T]$, one after the other. Looking at \Cref{lem:every edge twice} we obtain many constructions of 2-rainbow cycles:

\begin{figure}
        \centering
        \includegraphics[width=0.7\linewidth]{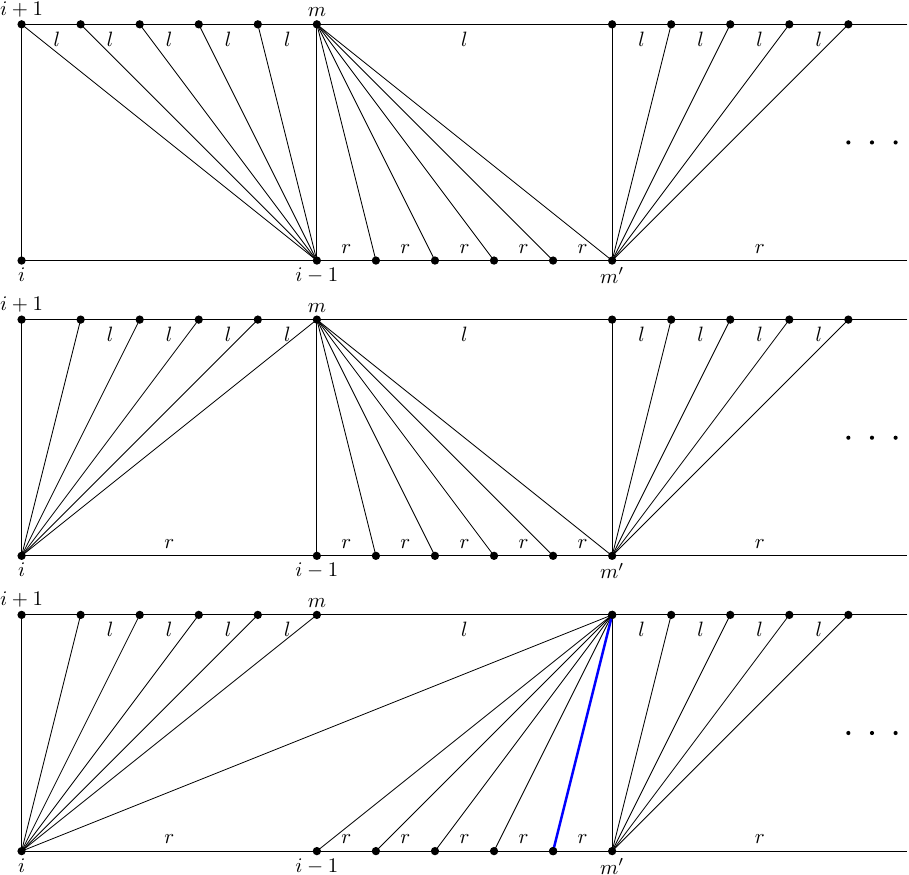}
        \caption{Illustration for the proof of \Cref{lem:rotation}. At the top, the initial zigzag triangulation $T$. In the middle, the triangulation after flipping the edges of $S_1^{[\:)}$. At the bottom, the triangulation after also flipping $\overline{S_2^{[\:)}}$. The letters of the word associated with the endpoint $i$ are written to the left of the associated edge.}
        \label{fig:rotation induction}
    \end{figure}

\begin{lemma}\label{lem:rotation}
    Let $T$ be a zigzag-triangulation with end point $i$ and $d$ a direction of rotation (clockwise or counter clockwise). Then there exists a permutation $\pi$ such that $T(\pi)$ equals $T$ rotated by one vertex in direction $d$ and hence
    \[
        [T] = \{T(\pi^k)\mid k\in\N\}.
    \]
    Furthermore, during the transformation, the word $\omega_i(T)$ is transformed blockwise from left to right such that the words associated to the intermediate triangulations differ from $\omega_i(T)$ in at most three positions which are in consecutive blocks. In particular, all intermediate triangulations are zigzag triangulations.
\end{lemma}
\begin{proof}
See \Cref{fig:rotation induction} for an illustration.
    Assign the edges of $T$ to the (ordered) sets $S_k$ in the following way. When constructing $\omega_i(T)$ we are walking across the dual path starting in $i$. When adding a letter, we are inside a triangle of $T$ and the letter we add depends which new vertex of the $n$-gon we see compared to the previous triangle. We associate this triangle to the letter. Given a block $b_j$ of $\omega_i(T)$, we identify the triangles in $T$ associated to the letters of this block. These edges form $S_j$ and are given the order they appear on the dual path when walking from $i$. Note that the sets $S_j$ are hence also indexed with respect to their order on the dual path starting from $i$. Note that $S_k$ and $S_{k+1}$ share an edge, i.e. the last edge of $S_k$ and the first edge of $S_{k+1}$ are equal. We set $S_k^{[\:]}$ to be $S_k$, $S_k^{(\:]}$ to be $S_k$ without its first edge, $S_k^{[\:)}$ to be $S_k$ without its last edge and $S_k^{(\:)}$ to be $S_k$ without its first and last edge. Furthermore, a set overset with a bar, for example $\overline{S_k^{(\:]}}$, represents the same set of edges with its ordering reversed. We claim that the following two orderings rotate $T$ by one vertex. One of them rotates clockwise and the other counter clockwise. If $b$ is even, then the orders are:
    \begin{align*}
        &S_1^{[\:)}\overline{S_2^{[\:)}}S_3^{[\:)}\overline{S_4^{[\:)}}\cdots S_{b-1}^{[\:)}\overline{S_b^{[\:]}}\\
        &\overline{S_1^{[\:]}}S_2^{(\:]}\overline{S_3^{[\:]}}S_4^{(\:]}\cdots \overline{S_{b-1}^{(\:]}}S_b^{(\:]}.
    \end{align*}
    If $b$ is odd, then the orders are:
    \begin{align*}
        &S_1^{[\:)}\overline{S_2^{[\:)}}S_3^{[\:)}\overline{S_4^{[\:)}}\cdots \overline{S_{b-1}^{[\:)}}S_b^{[\:]}\\
        &\overline{S_1^{[\:]}}S_2^{(\:]}\overline{S_3^{[\:]}}S_4^{(\:]}\cdots S_{b-1}^{(\:]}\overline{S_b^{(\:]}}.
    \end{align*}
    Note that in both cases the overline is supposed to alternate between the sets.

    First of all, note that the above really are permutations of the edges of $T$. We show the claimed properties under the assumption that $\pi = S_1^{[\:)}\overline{S_2^{[\:)}}\cdots S_{b-1}^{[\:)}\overline{S_b^{[\:]}}$. The proofs of the other cases are analogous. Without loss of generality, we will assume that $\omega_i(T) = r\ldots rl\ldots l\ldots$, i.e. that the first block is a block of $r$'s. We apply an induction on $b$. If $b=1$, then $T$ is a star and the result follows easily. If $b=2$, then, flipping the edges of $S_1^{[\:)}$ results in a sequence of zigzag triangulations with end point $i+1$ and words 
    \begin{align*}
        rrr\ldots rr&|ll\ldots lll\\
        lrr\ldots rr&|ll\ldots lll\\
        rlr\ldots rr&|ll\ldots lll\\
        &\vdots\\
        rrr\ldots rl&|ll\ldots lll.
    \end{align*}
    Next, flipping the edges $\overline{S_2^{[\:)}}$ transforms the word as follows:
    \begin{align*}
        rrr\ldots rl&|ll\ldots lll\\
        rrr\ldots rl&|ll\ldots llr\\
        rrr\ldots rl&|ll\ldots lrl\\
        &\vdots\\
        rrr\ldots rr&|ll\ldots lll.
    \end{align*}
    Note that all of these differ from the initial word in at most two places.
    
    Now assume that the result is shown for $b-2$ blocks. After performing the flips $S_1^{[\:)}\overline{S_2^{[\:)}}$, all flipped edges are in their desired position as dictated by the rotation. The last edge $e$ of $S_2^{[\:)}$ cuts the polygon, and we consider the zigzag triangulation obtained by restricting to the smaller polygon with the not yet flipped edges. Then, one vertex of $e$ is an end point of this restricted triangulation, and we see that the word associated to this end point consists of $b-2$ blocks, so we can apply induction to also rotate the remaining part of~$T$. Also note that, when performing flips in $S_k$, the edges of all sets $S_{k'}$ with~$k'\le k-2$ have been completely rotated, and hence the words of the intermediate triangulations agree with $\omega_i(T)$ in the associated positions. Therefore, it is not possible that letters from four different blocks disagree in any intermediate triangulation. Furthermore, considering that, when performing the flips from a set, we only have a single disagreement in that set, we see that there can generally be at most 3 disagreements between $\omega_i(T)$ and the intermediate words.
\end{proof}
    
\begin{lemma}\label{lem:2rainbow}
    Let $T$ be any zigzag triangulation without rotational symmetries and let $\pi$ be a permutation of its edges obtained by \Cref{lem:rotation}. Then $T\xrightarrow{\pi^{n}}T$ is a 2-rainbow cycle on $\A_n$.
\end{lemma}
\begin{proof}
    It is easy to see that we do not repeat any triangulation when rotating by one, as we flip every edge exactly once. The triangulations in different rotation steps differ by their endpoints and as there is no rotational symmetry the triangulations $T(\pi^ k)$, $k=1,\ldots,n$ are all different. Therefore $T\xrightarrow{\pi^{n}}T$ is a cycle. The edges introduced along its flips are exactly the edges in $T(\pi^k)$ for $k=1,\ldots,n$ because in $T(\pi^{k-1})\xrightarrow{\pi}T(\pi^k)$ every edge is flipped exactly once. By \Cref{lem:every edge twice}, the result follows.
\end{proof}

We now restrict our attention to a special class of zigzag triangulations we call \defn{imbalanced}. A zigzag triangulation is imbalanced if one (and therefore both) of its associated words are of the form $\ell r^m\ell$ or $l^m\bar\ell$, where $m=\lfloor\frac{n-4}{2}\rfloor$ and $\ell$ is a word of length $\lceil\frac{n-4}{2}\rceil$ which contains at least one $l$ and at least one $r$. Let $i$ be the endpoint of an imbalanced zigzag triangulation $T$ such that $\omega_i(T)=\ell r^m$ or $\omega_i(T)=\ell l^m$. Then we call $\ell$ the \defn{label} of $T$ and $i$ the \defn{pivot} of $T$. Let \defn{$\omega(T)=\omega_i(T)$} be the word associated to $T$ with respect to its pivot.
We need one more lemma:
\begin{lemma}\label{lem:labels}
    Let $m\in\N$ be given. There is a set of $k$ words on $\{r,l\}$ of length $m$ and
    \[
        k\le\frac{2^{m-1}}{\binom{m-1}{6}(2^6-1)}-2
    \]
    with the following properties:
    \begin{itemize}
        \item The pairwise differ in at least 7 positions.
        \item The difference of the number of words ending in $r$ versus the words ending in $l$ is at most 1.
        \item They differ in at least 7 positions from $r^m$ and $l^m$.
    \end{itemize}
    
\end{lemma}
\begin{proof}
    We pick words of length $m-1$ greedily and delete all words which differ in at most 6 positions from the set of chosable words until it becomes empty. For a given word, we have to delete $\binom{m-1}{6}(2^6-1)$ words (some of which might have been deleted before). Therefore, we can choose at least $k$ words if
    \[
        \binom{m-1}{6}(2^6-1)k\le 2^{m-1}\iff k\le \frac{2^{m-1}}{\binom{m-1}{6}(2^6-1)}.
    \]
    Then we can freely extend these words such that half ends with $r$ and the others with $l$ and maintain that they are pairwise different in at least 7 positions. If we choose $r^m$ and $l^m$ as the first words, then the others differ in at least 7 positions from these two as well.
\end{proof}

\begin{figure}
    \centering
    \scalebox{0.6}{

\def\drawlabels{}

% DEFCOORDS: Defines the coordinates used to draw triangulations
% 1: rotation
% 2: radius
% 3: N
\newcommand{\DEFCOORDS}[3]{
    \pgfmathsetmacro\stepangle{360/#3}
    \foreach\a in {1,...,#3}
    {
        \coordinate (\a) at (\a*\stepangle+#1*\stepangle-\stepangle:#2);
        \if\drawlabels y{\node at (\a*\stepangle+#1*\stepangle-\stepangle:#2+0.3) {\a};}\fi
    }
}

% TRI: Draws a triangulation with the specified parameters:
% 1: rotation
% 2: xshift
% 3: yshift
% 4: radius
% 5: parameters for the diagonals (like color for example)
% 6: N
% 7: list of diagonals as 2/4, ... 
\newcommand{\TRI}[7]{
    \begin{scope}[xshift=#2*#4*3cm, yshift=-#3*#4*3cm]
        \DEFCOORDS{#1}{#4}{#6}
        \foreach\x/\y in {#7}
        {
            \draw[#5, thick] (\x) -- (\y);
        }
        \pgfmathsetmacro\Nminusone{#6 -1}
        \foreach\a in {1,...,\Nminusone}
        {
            \pgfmathsetmacro\aplusone{\a+1}
            \draw (\a) -- (\aplusone);
            \fill (\a) circle (1pt);
        }
        \draw (1) -- (#6);
        \fill (#6) circle (1pt);
    \end{scope}
}

\def\radius{1}

\begin{tikzpicture}
    % n even:
    
    \TRI{0}{0}{0}{\radius}{very thick}{14}{1/8}
    \TRI{0}{0}{0}{\radius}{thick, blue}{14}{2/8, 2/7, 2/6, 2/5, 2/4}
    \TRI{0}{0}{0}{\radius}{thick, red}{14}{1/9, 9/14, 14/10, 10/13, 13/11}

    \draw[thick, ->] (1) +(5mm, 0) -- +(35mm,0) node [below, pos=0.5] {2-rainbow path} node [above, pos=0.5] {$\pi_1^{n-1}$};

    \TRI{-1}{2}{0}{\radius}{very thick}{14}{1/8}
    \TRI{-1}{2}{0}{\radius}{thick, blue}{14}{2/8, 2/7, 2/6, 2/5, 2/4}
    \TRI{-1}{2}{0}{\radius}{thick, red}{14}{1/9, 9/14, 14/10, 10/13, 13/11}

    \draw[thick, ->] (2) +(2mm,0) -- +(7mm,0);
    
    \TRI{-1}{3}{0}{\radius}{very thick}{14}{2/9}
    \TRI{-1}{3}{0}{\radius}{thick, blue}{14}{2/8, 2/7, 2/6, 2/5, 2/4}
    \TRI{-1}{3}{0}{\radius}{thick, red}{14}{1/9, 9/14, 14/10, 10/13, 13/11}

    \draw[thick, ->, blue] (2) +(5mm,0) -- +(20mm,0) node [above, pos=0.5] {$S\to T'_{2}$};

    \TRI{-1}{4.5}{0}{\radius}{very thick}{14}{2/9}
    \TRI{-1}{4.5}{0}{\radius}{thick, blue}{14}{9/3, 9/4, 4/8, 4/7, 4/6}
    \TRI{-1}{4.5}{0}{\radius}{thick, red}{14}{1/9, 9/14, 14/10, 10/13, 13/11}

    \draw[thick, ->, red] (2) +(5mm,0) -- +(20mm,0) node [above, pos=0.5] {$T'_1\to S$};
    
    \TRI{-1}{6}{0}{\radius}{very thick}{14}{2/9}
    \TRI{-1}{6}{0}{\radius}{thick, blue}{14}{9/3, 9/4, 4/8, 4/7, 4/6}
    \TRI{-1}{6}{0}{\radius}{thick, red}{14}{2/10, 2/11, 2/12, 2/13, 2/14}

    % and back

    \TRI{-1}{0}{1}{\radius}{very thick}{14}{2/9}
    \TRI{-1}{0}{1}{\radius}{thick, blue}{14}{9/3, 9/4, 4/8, 4/7, 4/6}
    \TRI{-1}{0}{1}{\radius}{thick, red}{14}{2/10, 2/11, 2/12, 2/13, 2/14}

    \draw[thick, ->] (2) +(5mm, 0) -- +(35mm,0) node [below, pos=0.5] {2-rainbow path} node [above, pos=0.5] {$\pi_{2}^{n-1}$};

    \TRI{0}{2}{1}{\radius}{very thick}{14}{2/9}
    \TRI{0}{2}{1}{\radius}{thick, blue}{14}{9/3, 9/4, 4/8, 4/7, 4/6}
    \TRI{0}{2}{1}{\radius}{thick, red}{14}{2/10, 2/11, 2/12, 2/13, 2/14}

    \draw[thick, ->, red] (1) +(5mm,0) -- +(20mm,0) node [above, pos=0.5] {$S\to T'_1$};

    \TRI{0}{3.5}{1}{\radius}{very thick}{14}{2/9}
    \TRI{0}{3.5}{1}{\radius}{thick, blue}{14}{9/3, 9/4, 4/8, 4/7, 4/6}
    \TRI{0}{3.5}{1}{\radius}{thick, red}{14}{1/9, 9/14, 14/10, 10/13, 13/12}
    
    \draw[thick, ->, blue] (1) +(5mm,0) -- +(20mm,0) node [above, pos=0.5] {$T_2\to S$};

    \TRI{0}{5}{1}{\radius}{very thick}{14}{2/9}
    \TRI{0}{5}{1}{\radius}{thick, blue}{14}{2/8, 2/7, 2/6, 2/5, 2/4}
    \TRI{0}{5}{1}{\radius}{thick, red}{14}{1/9, 9/14, 14/10, 10/13, 13/12}
    
    \draw[thick, ->] (1) +(2mm,0) -- +(7mm,0);

    \TRI{0}{6}{1}{\radius}{very thick}{14}{1/8}
    \TRI{0}{6}{1}{\radius}{thick, blue}{14}{2/8, 2/7, 2/6, 2/5, 2/4}
    \TRI{0}{6}{1}{\radius}{thick, red}{14}{1/9, 9/14, 14/10, 10/13, 13/12}

    % n odd

    \TRI{0}{0.5}{2.2}{\radius}{thick, blue}{15}{1/8, 1/7, 1/6, 1/5, 1/4, 1/3}
    \TRI{0}{0.5}{2.2}{\radius}{thick, red}{15}{1/9, 9/15, 15/10, 10/14, 14/11, 11/13}

    \draw[thick, ->] (1) +(5mm, 0) -- +(35mm,0) node [below, pos=0.5] {2-rainbow path} node [above, pos=0.5] {$\pi_1^{n-1}$};

    \TRI{-1}{2.5}{2.2}{\radius}{thick, blue}{15}{1/8, 1/7, 1/6, 1/5, 1/4, 1/3}
    \TRI{-1}{2.5}{2.2}{\radius}{thick, red}{15}{1/9, 9/15, 15/10, 10/14, 14/11, 11/13}

    \draw[thick, ->, blue] (2) +(5mm,0) -- +(20mm,0) node [above, pos=0.5] {$S\to T'_{2}$};

    \TRI{-1}{4}{2.2}{\radius}{thick, blue}{15}{2/9, 9/3, 9/4, 9/5, 5/8, 5/7}
    \TRI{-1}{4}{2.2}{\radius}{thick, red}{15}{1/9, 9/15, 15/10, 10/14, 14/11, 11/13}

    \draw[thick, ->, red] (2) +(5mm,0) -- +(20mm,0) node [above, pos=0.5] {$T'_1\to S$};

    \TRI{-1}{5.5}{2.2}{\radius}{thick, blue}{15}{2/9, 9/3, 9/4, 9/5, 5/8, 5/7}
    \TRI{-1}{5.5}{2.2}{\radius}{thick, red}{15}{2/10, 2/11, 2/12, 2/13, 2/14, 2/15}

    % and back

    \TRI{-1}{0.5}{3.2}{\radius}{thick, blue}{15}{2/9, 9/3, 9/4, 9/5, 5/8, 5/7}
    \TRI{-1}{0.5}{3.2}{\radius}{thick, red}{15}{2/10, 2/11, 2/12, 2/13, 2/14, 2/15}
    
    \draw[thick, ->] (2) +(5mm, 0) -- +(35mm,0) node [below, pos=0.5] {2-rainbow path} node [above, pos=0.5] {$\pi_2^{n-1}$};
    
    \TRI{0}{2.5}{3.2}{\radius}{thick, blue}{15}{2/9, 9/3, 9/4, 9/5, 5/8, 5/7}
    \TRI{0}{2.5}{3.2}{\radius}{thick, red}{15}{2/10, 2/11, 2/12, 2/13, 2/14, 2/15}
    
    \draw[thick, ->, red] (1) +(5mm,0) -- +(20mm,0) node [above, pos=0.5] {$S\to T'_1$};
    
    \TRI{0}{4}{3.2}{\radius}{thick, blue}{15}{2/9, 9/3, 9/4, 9/5, 5/8, 5/7}
    \TRI{0}{4}{3.2}{\radius}{thick, red}{15}{1/9, 9/15, 15/10, 10/14, 14/11, 11/13}
    
    \draw[thick, ->, blue] (2) +(5mm,0) -- +(20mm,0) node [above, pos=0.5] {$T'_2\to S$};
    
    \TRI{0}{5.5}{3.2}{\radius}{thick, blue}{15}{1/8, 1/7, 1/6, 1/5, 1/4, 1/3}
    \TRI{0}{5.5}{3.2}{\radius}{thick, red}{15}{1/9, 9/15, 15/10, 10/14, 14/11, 11/13}

    \draw[gray] (-\radius,-1.6*\radius*3 cm) -- (6*\radius*3 cm + \radius*1cm,-1.6*\radius*3 cm) node [above, pos=0.5] {$n$ even} node [below, pos=0.5] {$n$ odd};
    
\end{tikzpicture}
}
    \vspace{-9mm}
    \caption{An illustration for the proof of  \Cref{thm:associahedron}. Above, the case where $n$ is even. Below, the case where $n$ is odd.}
    \label{fig:construction imbalanced}
\end{figure}

We are ready to prove \Cref{thm:associahedron}:

\associahedron*
\begin{proof}
    Let $\ell_1,\ldots,\ell_N$ be words of length $m$ which pairwise differ in at least 7 positions and which also differ in at least 7 positions from $r^{\lfloor\frac{n-4}{2}\rfloor}$ and $l^{\lfloor\frac{n-4}{2}\rfloor}$, as given by \Cref{lem:labels}, with $m=\lceil\frac{n-4}{2}\rceil$. Assume that $\ell_i$ ends with $r$ if $i$ is odd, and ends in $l$ id $i$ is even. 

    First, assume that $n$ is even. A triangulation has $n-3$ diagonals, and a unique diagonal splitting the $n$-gon into two parts of the same size. This diagonal is called the center diagonal. We imagine it as drawn horizontally, between the vertices 1 and $\frac{n}{2}+1$, with 1 being on the right and $\frac{n}{2}-1$ on the left. Construct the triangulations $T_i$, $i=1,\ldots,N$ as follows: 
    \begin{itemize}
        \item If $i$ is odd, then insert all diagonals connecting 2 to the vertices $4,\ldots,\frac{n}{2}+1$ (remember that the vertices of the $n$-gon are labeled in counterclockwise order). This subtriangulation is a star with center at 2. Below the center diagonal insert diagonals so that the label of the resulting imbalanced zigzag triangulation has label $\ell_i$.
        \item If $i$ is even, then insert the diagonals connecting 1 to $n-1,\ldots,\frac{n}{2}+2$. This subtriangulation is a star with center at 1 below the center diagonal. Above the center diagonal, insert the diagonals such that the resulting imbalanced zigzag triangulation has label $\ell_i$.
    \end{itemize}
    The subtriangulation opposite to the star in $T_i$ (which gives its the label) is called $T_i'$.
    
    Now, let $\pi_i$ be the permutation given by \Cref{lem:rotation} to rotate $T_i$ counterclockwise from its pivot. Furthermore, \Cref{lem:adj to stars} guarantees the existence of a path flipping a star to $T'_i$ and back.

    The following is illustrated in \Cref{fig:construction imbalanced}.
    Starting from $T_1$, first apply $\pi_1^{n-1}$. This is a 2-rainbow path and results in the triangulation $T_1$ rotated by 1 vertex clockwise. Next, flip the center diagonal. Note that, because $\ell_1$ ends in $r$, this flip rotates the center diagonal by 1 vertex counterclockwise and the resulting triangulations is still zigzag. Next, transform the star above the center diagonal to $T'_2$ and then transform $T'_1$ to a star below the center diagonal. Note that, since the center of the star is at vertex 1 of the $n$-gon after the rotation, the subtriangulation $T'_2$ has an endpoint at 1. This implies that the labels match up. If $\ell_2$ would not end in $l$, then this transformation would not result in a zigzag triangulation. This results in the triangulation $T_2$.

    Next, rotate $T_2$ by applying $\pi_2^{n-1}$, transform the star to $T'_3$, then transform $T'_2$ to the star and finally flip the center diagonal. This procedure is repeated until we arrive at $T_N$. Because $N$ is even, we can use the same construction to flip back to $T_1$. In the figure, we flip to $T_2$ and then back, i.e. $N=2$.

    Now assume that $n$ is odd. In this case, there is no center diagonal. However, a very similar construction can be used, see \Cref{fig:construction imbalanced}. In this case, we construct $T_i$ as follows:
    \begin{itemize}
        \item If $i$ is odd, add the diagonals connecting 1 to $3,\ldots,\frac{n+11}{2}$, and then add the diagonals below this star so as to obtain the label $\ell_i$.
        \item If $i$ is even, add the diagonal connecting 1 to $n-2,\ldots,\frac{n-1}{2}$, and then add the diagonals below this star so as to obtain the label $\ell_i$.
    \end{itemize}
    The subtriangulation opposite to the star in $T_i$ (which gives its the label) is called $T_i'$.
    
    Now, let $\pi_i$ be the permutation given by \Cref{lem:rotation} to rotate $T_i$ counterclockwise from its pivot. Furthermore, \Cref{lem:adj to stars} guarantees the existence of a path flipping a star to $T'_i$ and back.

    Starting from $T_1$, apply $\pi_1^{n-1}$. This results in $T_1$ rotated by 1 vertex clockwise. Next, transform the star to $T'_2$ and then $T'_1$ to the star. This results in $T_2$. From $T_2$, apply $\pi_2^{n-1}$ which results in $T_2$ rotated by 1 vertex clockwise. Then transform the star to $T'_3$ and then $T'_2$ to the star to obtain $T_3$, and so on.

    In both cases of $n$ being even or odd, the above defined walk $C$ on the associahedron consists of $N$ 2-rainbow paths connected together by transforming the two subtriangulations one after the other, and possibly flipping the center diagonal. Note that this transformation connecting the 2-rainbow paths flips each diagonal of the end of a path exactly once. Hence, each diagonal of the $n$-gon appears in $C$ exactly $4N$ times. It thus only remains to show that $C$ is a cycle, i.e. that it does not repeat any triangulation. To show this, let $T$ be any triangulation in $C$. We will identify its position uniquely. 
    
    Assume first that $T$ is not zigzag, but contains a triangle $t$ not sharing a boundary edge with the $n$-gon. Then $T$ must be part of a transformation of a subtriangulation to or from a star. Therefore, the $n$-gon can be split into 2 pieces of the same size such that $t$ is fully contained in one of them and the other has exactly one vertex of degree 2. This vertex $i$ of degree 2 is the end point of the zigzag triangulations on the ends of this transformation. Consider the weak dual of $T$. When looking from $i$ we see a path of length at least $\lfloor\frac{n-4}{2}\rfloor$. Because we know that $T$ is a triangulation in the middle of a transformation involving the star part of an imbalanced zigzag triangulation, we can read of the label of the triangulation and hence recover the exact position of $T$ in $C$.

    Assume now that $T$ is zigzag. Then we can construct the two words $w_1$ and $w_2$ associated to $T$. Now match the labels $\ell_1,\ldots,\ell_N$ to the prefixes of $w_1,w_2$. The only way that both $w_1$ and $w_2$ can have a match is if they match consecutive labels $\ell_i,\ell_{i+1}$, and $T$ is exactly the triangulation obtained from $T_i$ after transforming the star to $T'_{i+1}$. Otherwise, exactly one must match because in each step, at least one label-defining set of diagonals is present. In both cases, we can identefy the labels and deduce the exact position of $T$ in $C$.
\end{proof}

\begin{remark}\label{rmk:smaller cycle}
    The lower bound on $n$ from the theorem can be reduced significantly if one is only interested in $r$-rainbow cycles for $r\in[2n+2]$. To achieve this, one can remove the star diagonals and consider more general zigzag triangulations. Then the corresponding 2-rainbow cycles can be glued together by flipping to an intermediary star and then back to a new zigzag triangulation. This way, more triangulations are available for the 2-rainbow cycles, reducing the lower bound on $n$. However, because there are only $n$ different stars, this fails for larger cycles and can only give a linear relation between $r$ and $n$. See \Cref{fig:rainbow big example} for an illustration of this construction.
\end{remark}

\begin{figure}[hbt]
        \newcommand{\defchords}[1]{
    \foreach\a in {1,...,10}{
        \coordinate (\a) at (-#1*36+\a*36-36:1);
        \draw[gray] (\a) -- (-#1*36+\a*36:1);
        \node at (\a*36-36:1.25) {\tiny\a};
    }
}

\def\xshift{2cm}
\def\yshift{3.5cm}

\newcommand{\Star}[4]{
\begin{scope}[xshift=#2*\xshift, yshift=#3*\yshift, scale=#4]
    \defchords{#1}
    \node[outer sep=2pt] (Star#1) at (0,0) {
        \begin{tikzpicture}
            \foreach\x/\y in {1/2,1/3,1/4,1/5,1/6,1/7,1/8,1/9}{
                \draw[thick] (\x) -- (\y);
            }
            \foreach\a in {1,...,10}{
                \fill (\a) circle (1pt);
            }
        \end{tikzpicture}
    };
\end{scope}
}

\newcommand{\Tone}[4]{
\begin{scope}[xshift=#2*\xshift, yshift=#3*\yshift, scale=#4]
    \defchords{#1}
    \node[outer sep=2pt] (Tone#1) at (0,0) {
        \begin{tikzpicture}
            \foreach\x/\y in {2/10, 2/9, 2/8, 8/3, 8/4, 4/7, 4/6}{
                \draw[thick] (\x) -- (\y);
            }
            \foreach\a in {1,...,10}{
                \fill (\a) circle (1pt);
            }
        \end{tikzpicture}
    };
\end{scope}
}

\newcommand{\Tri}[6]{
\begin{scope}[xshift=#2*\xshift, yshift=#3*\yshift, scale=#4]
    \defchords{#6}
    \node[outer sep=2pt] (#1) at (0,0) {
        \begin{tikzpicture}
            \foreach\x/\y in {#5}{
                \draw[thick] (\x) -- (\y);
            }
            \foreach\a in {1,...,10}{
                \fill (\a) circle (1pt);
            }
        \end{tikzpicture}
    };
\end{scope}
}

\newcommand{\Ttwo}[4]{
\begin{scope}[xshift=#2*\xshift, yshift=#3*\yshift, scale=#4]
    \defchords{#1}
    \node[outer sep=2pt] (Ttwo#1) at (0,0) {
        \begin{tikzpicture}
            \foreach\x/\y in {1/3, 3/9, 9/4, 4/8, 8/5, 8/4, 3/10, 5/7}{
                \draw[thick] (\x) -- (\y);
            }
            \foreach\a in {1,...,10}{
                \fill (\a) circle (1pt);
            }
        \end{tikzpicture}
    };
\end{scope}
}

\def\scale{0.6}
\def\innerscale{0.6}

\def\innerupper{-0.65}
\def\innerlower{-1.35}

\def\opacity{0.2}

\begin{center}
\scalebox{0.65}{
\begin{tikzpicture}
    \foreach\i in {0,...,9}{
        \pgfmathsetmacro\x{\i*(11/9)-2}
        \Tone{\i}{\x}{0}{\scale}
    }

    \Tri{ToneRot1}{3}{\innerupper}{\innerscale}{8/10, 1/8, 1/7, 7/2, 7/3, 3/6, 3/5}{2}
    \Tri{ToneRot2}{3}{\innerlower}{\innerscale}{8/10, 7/10, 1/7, 7/2, 7/3, 3/6, 3/5}{2}
    \Tri{ToneRot3}{4}{\innerlower}{\innerscale}{8/10, 7/10, 1/7, 7/2, 2/6, 3/6, 3/5}{2}
    \Tri{ToneRot4}{5}{\innerlower}{\innerscale}{8/10, 7/10, 1/7, 1/6, 2/6, 3/6, 3/5}{2}
    \Tri{ToneRot5}{5}{\innerupper}{\innerscale}{8/10, 7/10, 6/10, 1/6, 2/6, 3/6, 3/5}{2}
    \Tri{ToneRot6}{4}{\innerupper}{\innerscale}{8/10, 7/10, 6/10, 1/6, 2/6, 2/5, 3/5}{2}

    \Tri{ToneToStar1}{8}{\innerupper}{\innerscale}{1/9, 2/9, 2/8, 8/3, 8/4, 4/7, 4/6}{0}
    \Tri{ToneToStar2}{7}{\innerupper}{\innerscale}{1/9, 1/8, 2/8, 8/3, 8/4, 4/7, 4/6}{0}
    \Tri{ToneToStar3}{6}{\innerupper}{\innerscale}{1/9, 1/8, 1/3, 8/3, 8/4, 4/7, 4/6}{0}
    \Tri{ToneToStar4}{6}{\innerlower}{\innerscale}{1/9, 1/8, 1/3, 1/4, 8/4, 4/7, 4/6}{0}
    \Tri{ToneToStar5}{7}{\innerlower}{\innerscale}{1/9, 1/8, 1/3, 1/4, 1/7, 4/7, 4/6}{0}
    \Tri{ToneToStar6}{8}{\innerlower}{\innerscale}{1/9, 1/8, 1/3, 1/4, 1/7, 1/6, 4/6}{0}

    \Star{9}{9}{-1}{\scale}
    
    \foreach\i in {0,...,9}{
        \pgfmathsetmacro\x{(9-\i)*(11/9)-2}
        \Ttwo{\i}{\x}{-2}{\scale}
    }

    \Star{8}{-2}{\innerlower}{\innerscale}
    \Star{7}{-1}{\innerlower}{\innerscale}
    \Star{6}{0}{\innerlower}{\innerscale}
    \Star{5}{1}{\innerlower}{\innerscale}
    \Star{4}{2}{-1}{\innerscale}
    \Star{3}{1}{\innerupper}{\innerscale}
    \Star{2}{0}{\innerupper}{\innerscale}
    \Star{1}{-1}{\innerupper}{\innerscale}
    \Star{0}{-2}{\innerupper}{\innerscale}

    \draw[->] (Tone0) -- (Tone1) node[midway, above] {$\pi_1$};
    \draw[->] (Tone1) -- (Tone2) node[midway, above] {$\pi_1$};
    \draw[->] (Tone2) -- (Tone3) node[midway, above] {$\pi_1$};
    \draw[->] (Tone3) -- (Tone4) node[midway, above] {$\pi_1$};
    \draw[->] (Tone4) -- (Tone5) node[midway, above] {$\pi_1$};
    \draw[->] (Tone5) -- (Tone6) node[midway, above] {$\pi_1$};
    \draw[->] (Tone6) -- (Tone7) node[midway, above] {$\pi_1$};
    \draw[->] (Tone7) -- (Tone8) node[midway, above] {$\pi_1$};
    \draw[->] (Tone8) -- (Tone9) node[midway, above] {$\pi_1$};

    \draw[->] (Tone3) -- (ToneRot1);
    \draw[->] (ToneRot1) -- (ToneRot2);
    \draw[->] (ToneRot2) -- (ToneRot3);
    \draw[->] (ToneRot3) -- (ToneRot4);
    \draw[->] (ToneRot4) -- (ToneRot5);
    \draw[->] (ToneRot5) -- (ToneRot6);
    \draw[->] (ToneRot6) -- (Tone4);

    \draw[->] (Tone9) -- (ToneToStar1);
    \draw[->] (ToneToStar1) -- (ToneToStar2);
    \draw[->] (ToneToStar2) -- (ToneToStar3);
    \draw[->] (ToneToStar3) -- (ToneToStar4);
    \draw[->] (ToneToStar4) -- (ToneToStar5);
    \draw[->] (ToneToStar5) -- (ToneToStar6);
    \draw[->] (ToneToStar6) -- (Star9); 

    \draw[->] (Tone9) -- (Star9) node[midway, right] {$\pi'_1$};
    \draw[<-] (Star9) -- (Ttwo0) node[midway, right] {$\pi'_2$};

    \draw[->] (Ttwo0) -- (Ttwo1) node[midway, above] {$\pi_2$};
    \draw[->] (Ttwo1) -- (Ttwo2) node[midway, above] {$\pi_2$};
    \draw[->] (Ttwo2) -- (Ttwo3) node[midway, above] {$\pi_2$};
    \draw[->] (Ttwo3) -- (Ttwo4) node[midway, above] {$\pi_2$};
    \draw[->] (Ttwo4) -- (Ttwo5) node[midway, above] {$\pi_2$};
    \draw[->] (Ttwo5) -- (Ttwo6) node[midway, above] {$\pi_2$};
    \draw[->] (Ttwo6) -- (Ttwo7) node[midway, above] {$\pi_2$};
    \draw[->] (Ttwo7) -- (Ttwo8) node[midway, above] {$\pi_2$};
    \draw[->] (Ttwo8) -- (Ttwo9) node[midway, above] {$\pi_2$};

    \draw[->] (Ttwo9) -- (Star8) node[midway, right] {$\pi'_2$};
    
    \draw[->] (Star8) -- (Star7) node[midway, above] {$\pi$};
    \draw[->] (Star7) -- (Star6) node[midway, above] {$\pi$};
    \draw[->] (Star6) -- (Star5) node[midway, above] {$\pi$};
    \draw[->] (Star5) -- (Star4) node[midway, above] {$\pi$};
    \draw[->] (Star4) -- (Star3) node[midway, above] {$\pi$};
    \draw[->] (Star3) -- (Star2) node[midway, above] {$\pi$};
    \draw[->] (Star2) -- (Star1) node[midway, above] {$\pi$};
    \draw[->] (Star1) -- (Star0) node[midway, above] {$\pi$};

    \draw[<-] (Star0) -- (Tone0) node[midway, right] {$\pi'_1$};

    \fill[opacity=\opacity, color2] (Tone0.north west) -- (Tone9.north east) -- (Tone9.south east) -- (Tone0.south west) -- cycle;

    \fill[opacity=\opacity, pink] (Tone9.south east) -- (Tone9.south west) -- (Ttwo0.north west) -- (Ttwo0.north east) -- cycle;

    \fill[opacity=\opacity, color4] (Ttwo0.north east) -- (Ttwo9.north west) -- (Ttwo9.south west) -- (Ttwo0.south east) -- cycle;

    \fill[opacity=\opacity, color5] (Tone0.south west) -- (Tone0.south east) -- (Star3.north east) -- (Star4.north east) -- (Star4.south east) -- (Star5.south east) -- (Ttwo9.north east) -- (Ttwo9.north west) -- cycle;

    \fill[opacity=\opacity, color1] (ToneRot1.north west) -- (ToneRot5.north east) -- (ToneRot4.south east) -- (ToneRot2.south west) -- cycle;

    \fill[opacity=\opacity, color3] (ToneToStar1.north east) -- (ToneToStar3.north west) -- (ToneToStar4.south west) -- (ToneToStar6.south east) -- cycle;

\end{tikzpicture}}
\end{center}
        \vspace{-15mm}
    \caption{An illustration for the construction hinted at in \Cref{rmk:smaller cycle}. We assume that the zigzag triangulations $T_1, T_2$ are chosen so that that the following parts do not share triangulations. For bigger rainbow cycles, one would define $T_i$ for a larger range of $i$.}
    \vspace{5mm}
    \begin{minipage}{.5\textwidth}
        \begin{itemize}
            \item The first 2-rainbow path: \colorbox{color2!20}{$T_1\xrightarrow{\pi_1^{9}}T_1(\pi_1^{9})$}
            \item The connection: \colorbox{pink!20}{$T_1(\pi_1^9)\xrightarrow{\pi'_1}S_2\xleftarrow{\pi'_2}T_2$}
            \item The second 2-rainbow path: \colorbox{color4!20}{$T_2\xrightarrow{\pi_2^{9}}T_2(\pi_2^{9})$}
        \end{itemize}
    \end{minipage}
    \begin{minipage}{.55\textwidth}
        \begin{itemize}
            \item Visiting remaining stars:\\ \colorbox{color5!20}{$T_2(\pi_2^{9})\xrightarrow{\pi_2'}S_3\xrightarrow{\pi^{8}}S_1\xleftarrow{\pi_1'}T_1$}
            \item $\pi_1$ in detail: \colorbox{color1!20}{$T_1(\pi_1^{3})\xrightarrow{\pi_1}T_1(\pi_1^4)$}
            \item $\pi_1'$ in detail: \colorbox{color3!20}{$T_9\xrightarrow{\pi_1'}S_2$}
        \end{itemize}
    \end{minipage}
    \label{fig:rainbow big example}
\end{figure}

\pagebreak

\section{Permutation Gray Codes}\label{sec:permutation gray codes}

This section is dedicated to proving \Cref{thm:perm_rainbow}, \ref{thm:cyclicallyballanced} and \ref{thm: 2 and 3 rainbow for perm}.
The proof of \Cref{thm:perm_rainbow} will be done in three stages. First, we construct a balanced Gray code for $\Pall_n$ for odd $n\in \mathbb{N}$ in \Cref{sec:balancedOdd}. Next, a balanced Gray code for even $n$ is presented in \Cref{sec:balancedEven}. Finally, \Cref{thm:perm_rainbow} will follow by showing how to turn a balanced Gray code for $\Pall_n$ into a rainbow cycle on $\Pall_m$ with $m > n$ at the end of \Cref{sec:balancedEven}. The proof of \Cref{thm:cyclicallyballanced} is presented in \Cref{sec:cyclicallyBalanced} and the proof of \Cref{thm: 2 and 3 rainbow for perm} is presented in \Cref{sec:permutahedron}.

Following up on the preliminaries stated before the proof sketches, we now collect the remaining preliminaries for this section.
Subsequently, we sometimes write a permutation in \defn{cyclic notation} if we want to emphasize thinking of it as a bijective map ${\sigma:[n]\to [n]}$. For $\sigma\in S_{n}$ we say $\sigma$ \defn{acts on $S_{n}$} to define that
\[ \sigma(\pi) := \sigma(\pi_{1})\dotsc \sigma(\pi_{n}) \text{ for } \pi \in S_{n}. \]
On the other hand, for $\phi \in S_{n}$ we say $\phi$ \defn{acts on the indices} of $S_{n}$ to define
\[ \phi(\pi) := \pi_{\phi(1)}\dotsc \pi_{\phi(n)} \text{ for } \pi \in S_{n}. \]
As the terminology suggests both definitions give rise to a group action.

We say that \defn{a transposition $(i, j)$ contains $k$} if either $i=k$ or $j=k$. Furthermore, we frequently refer to the set of transpositions \defn{$T_{n}$} $:= \{(i, j) \mid i,j \in [n], i \not = j \}$ as the \defn{set of transpositions on $S_{n}$}. We sometimes say that \defn{$t$ is a transposition on $S_{n}$} to mean that~$t\in T_{n}$.

A \defn{partial} Gray code is a listing of the elements of a subset of $S_{n}$ where two consecutive elements still differ only by some transposition. A Gray code $\pi \xrightarrow{G} \pi'$ for $S_{n}$ is \defn{cyclic} if $\pi' \xrightarrow{t} \pi$ for some transposition $t$. To reiterate, we call a cyclic Gray code \defn{balanced} if each admissible transposition appears $r$ times within the cycle, where $r$ naturally depends on the admissible transpositions. We also call a transposition by itself balanced if it appears $r$ times within a cycle.
In \Cref{sec:balancedOdd} and \Cref{sec:balancedEven} the admissible transpositions are all transpositions on $S_{n}$, leading to $r = 2(n-2)!$, while in \Cref{sec:cyclicallyBalanced} only the cyclically adjacent ones are admissible, hence $r = (n-1)!$.

We will work with \defn{multisets} and quickly wish to establish a common ground on notation. A multiset $(A, \mu_{A})$ is a set $A$ together with a function $\mu_{A}:A \to \mathbb{N}^{+}$. For an element~$x\in A$ we call $\mu_{A}(x)$ the \defn{multiplicity} of $x$ and sometimes say that $x$ appears $\mu_{A}(x)$ times in $A$. We often leave the multiplicity function implicit and speak of the multiset $A$. The \defn{support} of a multiset $(A, \mu_{A})$ is defined to be the set~$A$. For clarity we list the definitions of operations we use on multisets $(A,\mu_{A})$ and $(B, \mu_{B})$:
\begin{itemize}
  \item $(A,\mu_{A}) \cupdot (B, \mu_{B}) := (A \cup B, \mu_{A\cup B})$ where $\mu_{A\cup B}(x) := \mu_{A}(x) + \mu_{B}(x)$.
  \item $(A,\mu_{A}) \setminus (B, \mu_{B}) := (A \setminus B, \mu_{A \setminus B})$ where $\mu_{A \setminus B}(x) :=  \max\{\mu_{A}(x) - \mu_{B}(x),0\}.$
  \item $|(A, \mu_{A})| := \sum_{a \in A}\mu_{A}(a)$.
\end{itemize}

If not specified otherwise, application of a function to a Gray code or multiset is always meant element-wise.

We now turn to proving \Cref{thm:perm_rainbow}.

\permrainbow*

As stated earlier, the proof is split into two parts. We first consider the case where $n$ is odd.

\subsection{A Balanced Gray Code for Permutations of Odd $[n]$} \label{sec:balancedOdd}
Let $n>1$ be odd. Our strategy for building a balanced Gray code for $\Pall_{n}$ is to reuse a balanced Gray code $G$ for $\Pall_{n-2}$ to produce an intermediate \emph{almost balanced Gray code} $L$ for $\Pall_{n-1}$, a notion we will define shortly. $L$ can then be utilized to generate a balanced Gray code for $\Pall_{n}$.

In order to count transpositions, we denote the number of times a transposition $t$ appears in a partial Gray code $G$ by $c_{t}^{G}$. This number is known as the \defn{transition count}. Note that for a cyclic Gray code our definition of transition count does not take the transposition closing the cycle into account. This is done in order to align it with concatenations of partial Gray codes.

Moving on, we want to measure how far transpositions in a Gray code $G$ are from being balanced. To this end, in order to measure the deviation from $2(n-2)!$, we define the multisets \defn{$c_{+}(G)$} by
\[  t \in c_{+}(G) \text{ with multiplicity } 2(n-2)! - c_t^G, \quad \text{ iff } c_t^G < 2(n-2)! \]
and \defn{$c_{-}(G)$} by
\[ t \in c_{-}(G) \text{ with multiplicity } c_t^G - 2(n-2)!, \quad \text{ iff } c_t^G > 2(n-2)! \]
where $t$ is a transposition on $S_{n}$. In case $t \not \in G$, we declare $c_{t}^{G} = 0$.
In other words,~$c_{+}(G)$ contains the transpositions that appear strictly less often than $2(n-2)!$ times in $G$ and~$c_{-}(G)$ respectively the ones appearing strictly more often, such that their deviation from $2(n-2)!$ corresponds to their multiplicity.
Using this notation we can rephrase the definition of an \defn{almost balanced Gray code}: A Gray code $G$ is almost balanced if it holds that $c_t^G = 2(n-2)!$
for all wide transpositions $t$ on $S_n$.

\begin{observation}
  If $\pi \xrightarrow{G} \pi'$ is a balanced Gray code for $\Pall_n$ then $c_{+}(G)$ is the singleton set containing the transposition $t$ with $\pi' \xrightarrow{t} \pi$, i.e. the transposition closing the cycle. Furthermore, there is by assumption no transposition appearing more often than $2(n-2)!$ in $G$ and, consequently, $c_{-}(G) = \emptyset$.
\end{observation}

We can see that for a balanced Gray code $G$ the multiset $c_{+}(G)$ trivially contains one more element then $c_{-}(G)$. This is also true in general as the next lemma will show.

\begin{lemma}\label{lem:T+T-relation}
  For any Gray code $G$ for $\Pall_n$ we have
  \[ |c_{-}(G)| + 1 = |c_{+}(G)| \]
  where $|\cdot|$ denotes the multiset's cardinality.
\end{lemma}
\begin{proof}
  During this proof, we treat the Gray code $G$ as a multiset of transpositions by simply forgetting the order in which they appear. By the definition of $c_{+}$ and $c_{-}$ each transposition on $S_{n}$ is contained in ${ (G \cupdot c_{+}(G)) \setminus c_{-}(G) }$ and has multiplicity $2(n-2)!$. As a result,
  \[ n! = |(G \cupdot c_{+}(G)) \setminus c_{-}(G)| = |G| + |c_{+}(G)| - |c_{-}(G)|. \]
  Furthermore, we know that $|G| = n! - 1$. Combining these equations yields
  \begin{align*}
         & n! - 1 + |c_{+}(G)| - |c_{-}(G)| = n!\\
    \iff \quad & |c_{-}(G)| + 1 = |c_{+}(G)|
  \end{align*}
  which completes the proof.
\end{proof}

In what follows, there will occasionally be the need to extend a permutation of $[n]$ to the set of transpositions on $S_{n}$. The following lemma covers this construction.

\begin{lemma}\label{lem:transbij}
  Let $\sigma:[n] \to [n]$ be a permutation. Then $\sigma$ extends to a bijection on the set of transpositions on $S_{n}$ by setting
  \[ \sigma((i, j)) := (\sigma(i), \sigma(j)) \]
  for a transposition $(i, j)$.
\end{lemma}
\begin{proof}
  Let $\sigma:[n] \to [n]$ be a permutation. Naturally, $\sigma$ extends to a bijection $\bar{\sigma}$ on~$[n] \times [n]$ if we define $\bar{\sigma}(i, j) := (\sigma(i), \sigma(j))$ for $i, j \in [n]$. Now consider the equivalence relation on~$[n] \times [n]$ generated by $(i, j) \sim (j, i)$.
  The map $\bar{\sigma}$ descends to a bijection $\hat{\sigma}:\nicefrac{[n] \times [n]}{\sim} \to \nicefrac{[n] \times [n]}{\sim}$ because \[\bar{\sigma}(i, j) = (\sigma(i), \sigma(j)) \sim  (\sigma(j), \sigma(i)) = \bar{\sigma}(j,i).\]
  Since \[(\nicefrac{[n] \times [n]}{\sim}) \setminus (\nicefrac{\{(i,i)|i\in [n]\}}{\sim}) \] is in one-to-one correspondence with the set of transpositions, the result follows.
\end{proof}

Because of \Cref{lem:transbij} we establish the convention that for a permutation $\sigma$ and a transposition $(i, j)$, the application \defn{$\sigma((i, j))$} will always be defined as \[(\sigma(i), \sigma(j)).\]

Next, we state two lemmas which we will use later when constructing the desired codes.

\begin{lemma}\label{lem:cycleapplication}
  Let $G$ be a partial Gray code for $\Pall_n$ and let further $w_{l,h} := (l \dotsc h) \in S_{n}$ act on $S_{n}$. Then $w_{l,h}(G)$ is a partial Gray code for $\Pall_n$. Furthermore, if no transposition in $G$ contains $k$, then no transposition in $\sigma_{l,h}(G)$ contains $\sigma_{l,h}(k)$.
\end{lemma}
\begin{proof}
  Let $\pi \in S_n$ be an arbitrary starting permutation for $w_{l,h}(G)$ and let $P$ denote the set of permutations listed this way. Note that $G$ when starting at $w_{l,h}^{-1}(\pi)$ lists the permutations $w_{l,h}^{-1}(P)$. By assumption all permutations in $w_{l,h}^{-1}(P)$ are distinct. Since~$w_{l,h}^{-1}$ is bijective, it follows that all permutations in $P$ are distinct. Hence, $w_{l,h}(G)$ is a partial Gray code.
  Now suppose that no transposition in $G$ contains $k$ and let $(i, j) \in G$. Since we have $\sigma_{l,h}((i, j)) = (\sigma_{l,h}(i), \sigma_{l,h}(j))$, we immediately see that no transposition in $\sigma_{l,h}(G)$ contains $\sigma_{l,h}(k)$.
\end{proof}

\begin{lemma}\label{lem:balancedCodeH}
  Let $12 \dotsc (n-1)n \xrightarrow{G} 12 \dotsc n(n-1)$ be a balanced Gray code for $\Pall_n$ with the property $c_{+}(G) = \{(n-1, n)\}$. Then there exists a balanced Gray code $H$ \[12 \dotsc (n-1)n \xrightarrow{H} 21 \dotsc (n-1)n\]
  with $c_{+}(H) = \{(1, 2)\}$ for $\Pall_n$.
\end{lemma}
\begin{proof}
  Define $\sigma := (1\dotsc n)$ acting on $S_{n}$ and $\phi := (1 \dotsc n)$ acting on the indices of $S_{n}$. As a result, 
  \[ \phi^{n-2}(12 \dotsc (n-1)n) \xrightarrow{G} \phi^{n-2}(12 \dotsc n(n-1)) \]
  is a balanced Gray code for $\Pall_n$. Using \Cref{lem:cycleapplication}, it is clear that 
  \[ \sigma^2(\phi^{n-2}(12 \dotsc (n-1)n)) \xrightarrow{\sigma^2(G)} \sigma^2(\phi^{n-2}(12 \dotsc n(n-1))) \]
  is also a balanced Gray code for $\Pall_n$. In particular, we have that 
  \begin{align*}
    &\sigma^{2}(\phi^{n-2}(12 \dotsc (n-1)n)) = 12 \dotsc (n-1)n, \\
    &\sigma^{2}(\phi^{n-2}(12 \dotsc (n-1)n)) = 21 \dotsc (n-1)n
  \end{align*}
  and using bijectivity of $\sigma$
  \[c_{+}(\sigma^2(G)) = \sigma^{2}(\{(n-1, n)\}) = \{(1, 2)\}.\] 
  Hence, setting $H := \sigma^2(G)$ proves the claim.
\end{proof}

We now present the almost balanced Gray code our construction relies on.

\begin{proposition}\label{prop:almostbalancedCode}
  Let $n\in \mathbb{N}$ be odd and let $12\dotsc (n-1)n \xrightarrow{G} 12\dotsc n(n-1)$ be a balanced Gray code for $\Pall_n$ such that $(n-1, n)$ is the transposition closing the cycle. Then there exists an almost balanced Gray code \[ 12\dotsc (n+1) \xrightarrow{L} 23\dotsc (n+1)1\]
  for $\Pall_{n+1}$.
\end{proposition}

%
% Creation of almost balanced Gray Code for S_4
%

\begin{figure}[h!] 
\centering
\begin{minipage}[t]{.25\textwidth}
\centering
% Gray Code G
\begin{tikzpicture}[rotate=270]
    \newcommand{\n}{3}                          %size of the permutations
    \newcommand{\width}{3}                       
    \newcommand{\height}{1.5}                       

    \pgfmathsetmacro\N{\n!}
    \pgfmathsetmacro\stepright{\width / \N}
    \pgfmathsetmacro\stepup{\height / \n}

%draw colors
\foreach\i/\a/\b/\c in                       %this needs to be expanded if \n is bigger
{1/1/2/3,
2/3/2/1,
3/3/1/2,
4/2/1/3,
5/2/3/1,
6/1/3/2
}
{
    \newcommand{\sethgt}[1]
    {
        \if##1\a\def\hgt{1}\fi
        \if##1\b\def\hgt{2}\fi
        \if##1\c\def\hgt{3}\fi
    }
        \foreach\y in {1,...,\n}
        {
            \sethgt{\y} %sets \hgt
            \setclr{\y} %sets \clr
            \fill[\clr] (\i * \stepright, \hgt*\stepup) -- (\i * \stepright + \stepright, \hgt*\stepup) -- (\i * \stepright + \stepright, \hgt*\stepup + \stepup) -- (\i * \stepright, \hgt*\stepup + \stepup) -- cycle;
        }
}

%grid
    \foreach\y in {0,...,\n}
    {
        \draw (\stepright, \y * \stepup + \stepup) -- (\N * \stepright + \stepright, \y * \stepup + \stepup);
    }
    \foreach\x in {0,...,\N}
    {
        \draw (\x * \stepright + \stepright, \stepup) -- (\x * \stepright + \stepright, \n * \stepup + \stepup);
    }

\end{tikzpicture}

$G$
\end{minipage}%
% Gray Code H
\begin{minipage}[t]{.25\textwidth}
\centering
% Gray Code H
\begin{tikzpicture}[rotate=270]
    \newcommand{\n}{3}                          %size of the permutations
    \newcommand{\width}{3}                       
    \newcommand{\height}{1.5}                       

    \pgfmathsetmacro\N{\n!}
    \pgfmathsetmacro\stepright{\width / \N}
    \pgfmathsetmacro\stepup{\height / \n}

%draw colors
\foreach\i/\a/\b/\c in                       %this needs to be expanded if \n is bigger
{1/1/2/3,
2/1/3/2,
3/3/1/2,
4/3/2/1,
5/2/3/1,
6/2/1/3
}
{
    \newcommand{\sethgt}[1]
    {
        \if##1\a\def\hgt{1}\fi
        \if##1\b\def\hgt{2}\fi
        \if##1\c\def\hgt{3}\fi
    }
        \foreach\y in {1,...,\n}
        {
            \sethgt{\y} %sets \hgt
            \setclr{\y} %sets \clr
            \fill[\clr] (\i * \stepright, \hgt*\stepup) -- (\i * \stepright + \stepright, \hgt*\stepup) -- (\i * \stepright + \stepright, \hgt*\stepup + \stepup) -- (\i * \stepright, \hgt*\stepup + \stepup) -- cycle;
        }
}

%grid
    \foreach\y in {0,...,\n}
    {
        \draw (\stepright, \y * \stepup + \stepup) -- (\N * \stepright + \stepright, \y * \stepup + \stepup);
    }
    \foreach\x in {0,...,\N}
    {
        \draw (\x * \stepright + \stepright, \stepup) -- (\x * \stepright + \stepright, \n * \stepup + \stepup);
    }

\end{tikzpicture}

$H$
\end{minipage}%
% Almost balanced Gray Code for S_4
\begin{minipage}[t]{.25\textwidth}
\centering
% Gray Code L
\begin{tikzpicture}[rotate=270]
    \newcommand{\n}{4}                          %size of the permutations
    \newcommand{\width}{10}                       
    \newcommand{\height}{1.6}                       

    \pgfmathsetmacro\N{\n!}
    \pgfmathsetmacro\stepright{\width / \N}
    \pgfmathsetmacro\stepup{\height / \n}

%draw colors
\foreach\i/\a/\b/\c/\d in                       %this needs to be expanded if \n is bigger
{1/1/2/3/4,
2/1/3/2/4,
3/3/1/2/4,
4/3/2/1/4,
5/2/3/1/4,
6/2/1/3/4,
7/2/1/4/3,
8/4/1/2/3,
9/1/4/2/3,
10/2/4/1/3,
11/4/2/1/3,
12/1/2/4/3,
13/1/3/4/2,
14/4/3/1/2,
15/4/1/3/2,
16/3/1/4/2,
17/3/4/1/2,
18/1/4/3/2,
19/2/4/3/1,
20/4/2/3/1,
21/4/3/2/1,
22/3/4/2/1,
23/3/2/4/1,
24/2/3/4/1
}
{
    \newcommand{\sethgt}[1]
    {
        \if##1\a\def\hgt{1}\fi
        \if##1\b\def\hgt{2}\fi
        \if##1\c\def\hgt{3}\fi
        \if##1\d\def\hgt{4}\fi
    }
        \foreach\y in {1,...,\n}
        {
            \sethgt{\y} %sets \hgt
            \setclr{\y} %sets \clr
            \fill[\clr] (\i * \stepright, \hgt*\stepup) -- (\i * \stepright + \stepright, \hgt*\stepup) -- (\i * \stepright + \stepright, \hgt*\stepup + \stepup) -- (\i * \stepright, \hgt*\stepup + \stepup) -- cycle;
        }
}

%grid
    \foreach\y in {0,...,\n}
    {
        \draw (\stepright, \y * \stepup + \stepup) -- (\N * \stepright + \stepright, \y * \stepup + \stepup);
    }
    \foreach\x in {0,...,\N}
    {
        \draw (\x * \stepright + \stepright, \stepup) -- (\x * \stepright + \stepright, \n * \stepup + \stepup);
    }

\end{tikzpicture}

$L$
\end{minipage}%

\caption{Visualization of Gray codes used in the proof of Proposition \ref{prop:almostbalancedCode} for $n=3$ ($1=$~purple, $2=$~blue, $3=$~green and $4=$~yellow). $L$ shows the almost balanced Gray code for~$\Pall_4$.}
\label{fig:almostbalancedL}
\end{figure}

\begin{proof}
  Taking the balanced Gray code $G$ as a starting point and applying \Cref{lem:balancedCodeH}, we obtain a second balanced Gray code
  
  \[ 12\dotsc (n-1)n \xrightarrow{H} 21\dotsc (n-1)n \]
  
  with $c_{+}(H) = \{(1, 2)\}$ and $c_{-}(H) = \emptyset$ for $\Pall_n$. Define $\omega_{k} := (k \dotsc n+1)$ acting on $S_{n+1}$. The diagram displayed in \Cref{fig:almostbalancedCode} shows the Gray code $L$ which we want to show to be almost balanced.

  \begin{figure}[h!]
    \centering
    \[
      \begin{tikzcd}[row sep=large]
        12\dotsc n (n+1) \arrow[rr, "H"]     &        & 21\dotsc n (n+1) \arrow[lld, "{(n, n+1)}"']   \\
        21\dotsc (n+1)n \arrow[rr, "\overleftarrow{\omega_{n}(H)}"]   &        & 12\dotsc (n+1)n \arrow[lld, "{(n-1, n)}"']  \\
        12\dotsc n(n+1)(n-1) \arrow[rrd, phantom, "\vdots" description]                       &  & \\
        & & 124 \dotsc n(n+1)3 \arrow[lld, "{(2, 3)}"'] \\
        134\dotsc n(n+1)2 \arrow[rr, "\omega_{2}(G)"]      &        & 134\dotsc (n+1)n2 \arrow[lld, "{(1, 2)}"']  \\
        234\dotsc (n+1)n1 \arrow[rr, "\overleftarrow{\omega_{1}(G)}"] &        & 234\dotsc n(n+1)1
      \end{tikzcd}
    \]
    \caption{Almost balanced Gray code for $\Pall_{n+1}$ with $n$ odd.}
    \label{fig:almostbalancedCode}
  \end{figure}

  Note that $n+1$ is not contained in any transposition found in $H$ or $G$. As a result of \Cref{lem:cycleapplication}, the codes $\omega_l(H)$ and $\omega_l(G)$ then do not change $\omega_l(n+1) = l$. Hence, for permutations ending with an even number $l$ for $4 \leq l \leq n+1$ the partial Gray code 
  \[ \omega_l(12\dotsc(n+1)) \xrightarrow{\omega_l(H)} \omega_l(21\dots(n+1)) \] 
  is used. Permutations ending with an odd number $k$ for $3 \leq k \leq n$ are listed by
  \[ \omega_k(21\dotsc(n+1)) \xrightarrow{\overleftarrow{\omega_k(H)}} \omega_k(12\dots(n+1)),\] 
 the permutations ending with $2$ are listed by 
 \[ \omega_2(1 \dotsc (n-1)n(n+1)) \xrightarrow{\omega_2(G)} \omega_2(1 \dots n(n-1)(n+1)) \]
 and the once ending with $1$ by 
 \[ \omega_1(1 \dots n(n-1)(n+1)) \xrightarrow{\overleftarrow{\omega_1(G)}} \omega_1(1 \dotsc (n-1)n(n+1)). \]
 
 For $k\in\{2,\dotsc,n+1\}$ the code listing the permutations ending with $k$ is connected to the code listing the permutations ending with $k-1$ by the transposition $(k-1, k)$. Consequently, $L$ is a Gray code for $\Pall_{n+1}$. We now prove that $L$ is almost balanced.

  Let $k\in \{3,\dotsc,n+1\}$. Because $H$ is balanced such that $c_{+}(H) = \{(1, 2)\}$ and $\omega_{k}$ is bijective with $\omega_{k}((1, 2)) = (1, 2)$, we further have that all transpositions on $S_{n+1}$ not containing $k$ appear $2(n-2)!$ times in $\omega_{k}(H)$ except for $(1, 2)$.

  Now let $k\in \{1,2\}$. By the same argument every transposition not containing $k$ is contained~$2(n-2)!$ times in $\omega_{k}(G)$ except for $(n, n+1)$. Hence, for $(i, j)$ a wide transposition on $S_{n+1}$, we have established that $(i, j)$ is contained
  \begin{itemize}
    \item $2(n-2)!$ times in $\omega_{k}(H)$ for $k \in \{3,\dotsc,n+1\} \setminus \{i, j\}$ and
    \item $2(n-2)!$ times in $\omega_{k}(G)$ for $k \in \{1,2\} \setminus \{i, j\}$.
  \end{itemize}
  In total, $(i, j)$ is contained $(n-1) \cdot 2(n-2)! = 2(n-1)!$ in $L$. As a result, $L$ is almost balanced.
\end{proof}

Given the almost balanced Gray code from the above proposition, we can now generate a balanced Gray code from it.

\begin{proposition}\label{prop:balancedCode}
  Let $n\in \mathbb{N} \setminus \{0\}$ be even and let $12\dotsc n \xrightarrow{L} 23\dotsc n1$ be an almost balanced Gray code for $S_{n}$. Then there exists a balanced Gray code
  \[12\dotsc n(n+1) \xrightarrow{G} 12\dotsc (n+1)n\]
  for $\Pall_{n+1}$ with the closing transposition being $(n, n+1)$.
\end{proposition}
\begin{proof}
  The idea for the proof is to apply the rotation $\sigma := (1 \dotsc n+1) \in S_{n+1}$, acting on $S_{n+1}$, to the almost balanced Gray code $L$ repeatedly. We will show that this process yields a balanced Gray code for $\Pall_{n+1}$.
  The diagram in \Cref{fig:balancedCode} shows the cyclic Gray code, call it $G$, for $\Pall_{n+1}$ which results from this approach.

  \begin{figure}[h]
    \centering
    \[
      \begin{tikzcd}
        12\dotsc n(n+1) \arrow[r, phantom] \arrow[rr, "L"] & {}     & 23\dotsc n1(n+1) \arrow[lld, "{(1, n+1)}"']   \\
        23\dotsc n(n+1)1 \arrow[rr, "\sigma(L)"]                       &        & 34\dotsc (n+1)21 \arrow[lld, "{(1, 2)}"']     \\
        34\dotsc (n+1)12 \arrow[rrd, phantom, "\vdots" description] &  & \\
        & & (n+1)1\dotsc n(n-1) \arrow[lld, "{(n-1, n)}"] \\
        (n+1)1\dotsc (n-1)n \arrow[rr, "\sigma^{n}(L)"']   &        & 12\dotsc (n+1)n
      \end{tikzcd}
    \]
    \caption{Balanced Gray code for $\Pall_n$ with $n$ odd.}
    \label{fig:balancedCode}
  \end{figure}

The permutations of $S_{n+1}$ ending with $k$ are listed by the partial Gray code 
  \[ \sigma^{k}(12\dotsc n(n+1)) \xrightarrow{\sigma^{k}(L)} \sigma^k(23\dotsc n1(n+1)). \] 
  Furthermore, the code $\sigma^{k}(L)$ is connected to $\sigma^{k+1}(L)$ by the transposition ${(\sigma^{k}(1), \sigma^{k}(n+1))}$.
  We first show that the wide transpositions without $(1, n+1)$ are balanced in $G$ and then turn to the cyclically adjacent ones.

  For the wide transpositions without $(1, n+1)$ the argument is similar to the one used in \Cref{prop:almostbalancedCode}. We observe using \Cref{lem:cycleapplication} that $\sigma^{k}(L)$ does not contain a transposition containing $k$ because no transposition in $L$ contains $n+1$. Furthermore, all wide transpositions on $S_{n+1}$ not containing $n+1$ appear $2(n-2)!$ times in $L$. It follows that all wide transpositions on $S_{n+1}$ not containing $k$ appear~$2(n-2)!$ times in $\sigma^{k}(L)$. As a result, a wide transposition $(i, j)$ on $S_{n+1}$, except $(1, n +1)$, is not contained in $\sigma^{i}(L)$ and $\sigma^{j}(L)$, while it appears $2(n-2)!$ times otherwise. Consequently, $(i, j)$ appears $2(n-1)!$ times in $G$ and is, hence, balanced.

  Until the end of this proof we treat the codes $G$ and $L$ as multisets of transpositions by forgetting the order. Denote by \[Z := \{\sigma^{k}((1, n+1)) \mid k \in \{0,\dotsc,n-1\}\] the transpositions connecting the partial Gray codes $\sigma^{k}(L)$ to $\sigma^{k+1}(L)$ as seen in the diagram. In order to prove the cyclically adjacent transpositions to be balanced, we start by noting that
  \begin{equation*}
      G \setminus Z = \bigcupdot_{k = 0}^{n} \sigma^{k}(L).
  \end{equation*}
  Furthermore, $\sigma$ and $\sigma^{-1}$ both preserve cyclically adjacent transpositions. We conclude that
  \begin{align*}
    & \{t \in G \mid t \text{ is cyclically adjacent}\} \setminus Z \\
    = \quad & \bigcupdot_{k = 0}^{n}\sigma^{k}(\{t \in L \mid t \text{ is adjacent}\}) \\
    = \quad & \bigcupdot_{\substack{t \in L\\ t \text{ is adjacent}}} \{\sigma^{k}(t) \mid k \in \{0,\dotsc,n\} \\
    = \quad & \bigcupdot_{\substack{t \in L\\ t \text{ is adjacent}}} \{\sigma^{k}((1, 2)) \mid k \in \{0,\dotsc,n\}\}.
  \end{align*}
  Hence, the cardinality of $\{t \in L \mid t \text{ is adjacent}\}$ determines the number of times a cyclically adjacent transposition on $S_{n+1}$ appears in $G \setminus Z $.
  To calculate this cardinality, we observe that for a transposition~$(i, j)$ in $L$ the transition count is given by
  \[ c_{(i, j)}^{L} = 2(n-2)! - \mu_{c_{+}(L)}((i, j)) + \mu_{c_{-}(L)}((i, j)) \]
  utilizing the definitions of $c_{-}$ and $c_{+}$. Denote by $N$ the set of adjacent transpositions on~$S_{n}$. It then follows that
  \begin{align*}
    & | \{t \in L \mid t \text{ is adjacent}\} | \\
    = \quad & \sum_{t \in N} c_{t}^{L} \\
    = \quad & \sum_{t \in N}2(n-2)! - \sum_{t \in N}\mu_{c_{+}(L)}(t) + \sum_{t \in N}\mu_{c_{-}(L)}(t) \\
    = \quad & \left(\sum_{i=1}^{n-1}2(n-2)!\right) - |c_{+}(L)| + |c_{-}(L)|
  \end{align*}
  The last line is implied by the assumption of $L$ being almost balanced, so that $c_{+}(L)$ and~$c_{-}(L)$ cannot contain a wide transposition. We now apply \Cref{lem:T+T-relation} to obtain
  \[ \left(\sum_{i=1}^{n-1}2(n-2)!\right) - |c_{+}(L)| + |c_{-}(L)| = (n-1) \cdot 2(n-2)! - 1 = 2(n-1)! - 1.\]
  In total this yields that every cyclically adjacent transposition on $S_{n+1}$ appears $2(n-1)! - 1$ times in~$G \setminus Z $. Hence, every cyclically adjacent transposition on $S_{n+1}$, except $(n, n+1)$, is balanced in $G$ because $(n, n+1)$ is the only cyclically adjacent transposition not included in $Z$. Observing, however, that $(n, n+1)$ is the transposition closing the cyclic Gray code $G$ implies that $G$ is, indeed, a balanced Gray code for $\Pall_{n+1}$.
\end{proof}

%
% Creation of balanced Gray Code for S_5
%

\begin{figure}[h!]
\begin{minipage}[t]{.19\textwidth}
\centering
% Almost balanced Gray Code L
\begin{tikzpicture}[rotate=270]
    \newcommand{\n}{5}                          %size of the permutations
    \newcommand{\width}{10}                       
    \newcommand{\height}{2}                       

    \pgfmathsetmacro\N{24}
    \pgfmathsetmacro\stepright{\width / \N}
    \pgfmathsetmacro\stepup{\height / \n}

%draw colors
\foreach\i/\a/\b/\c/\d/\e in                       %this needs to be expanded if \n is bigger
{1/1/2/3/4/5,
2/1/3/2/4/5,
3/3/1/2/4/5,
4/3/2/1/4/5,
5/2/3/1/4/5,
6/2/1/3/4/5,
7/2/1/4/3/5,
8/4/1/2/3/5,
9/1/4/2/3/5,
10/2/4/1/3/5,
11/4/2/1/3/5,
12/1/2/4/3/5,
13/1/3/4/2/5,
14/4/3/1/2/5,
15/4/1/3/2/5,
16/3/1/4/2/5,
17/3/4/1/2/5,
18/1/4/3/2/5,
19/2/4/3/1/5,
20/4/2/3/1/5,
21/4/3/2/1/5,
22/3/4/2/1/5,
23/3/2/4/1/5,
24/2/3/4/1/5
}
{
    \newcommand{\sethgt}[1]
    {
        \if##1\a\def\hgt{1}\fi
        \if##1\b\def\hgt{2}\fi
        \if##1\c\def\hgt{3}\fi
        \if##1\d\def\hgt{4}\fi
        \if##1\e\def\hgt{5}\fi
    }
        \foreach\y in {1,...,\n}
        {
            \sethgt{\y} %sets \hgt
            \setclr{\y} %sets \clr
            \fill[\clr] (\i * \stepright, \hgt * \stepup) -- (\i * \stepright + \stepright, \hgt * \stepup) -- (\i * \stepright + \stepright, \hgt * \stepup + \stepup) -- (\i * \stepright, \hgt * \stepup + \stepup) -- cycle;
        }
}

%grid
    \foreach\y in {0,...,\n}
    {
        \draw (\stepright, \y * \stepup + \stepup) -- (\N * \stepright + \stepright, \y * \stepup + \stepup);
    }
    \foreach\x in {0,...,\N}
    {
        \draw (\x * \stepright + \stepright, \stepup) -- (\x * \stepright + \stepright, \n * \stepup + \stepup);
    }   
    
\end{tikzpicture}

$L$
\end{minipage}%
% Almost balanced Gray Code \sigma(L)
\hfill
\begin{minipage}[t]{.19\textwidth}
\centering
% Almost balanced Gray Code \sigma(L)
\begin{tikzpicture}[rotate=270]
    \newcommand{\n}{5}                          %size of the permutations
    \newcommand{\width}{10}                       
    \newcommand{\height}{2}                       

    \pgfmathsetmacro\N{24}
    \pgfmathsetmacro\stepright{\width / \N}
    \pgfmathsetmacro\stepup{\height / \n}

%draw colors
\foreach\i/\a/\b/\c/\d/\e in                       %this needs to be expanded if \n is bigger
{1/2/3/4/5/1,
2/2/4/3/5/1,
3/4/2/3/5/1,
4/4/3/2/5/1,
5/3/4/2/5/1,
6/3/2/4/5/1,
7/3/2/5/4/1,
8/5/2/3/4/1,
9/2/5/3/4/1,
10/3/5/2/4/1,
11/5/3/2/4/1,
12/2/3/5/4/1,
13/2/4/5/3/1,
14/5/4/2/3/1,
15/5/2/4/3/1,
16/4/2/5/3/1,
17/4/5/2/3/1,
18/2/5/4/3/1,
19/3/5/4/2/1,
20/5/3/4/2/1,
21/5/4/3/2/1,
22/4/5/3/2/1,
23/4/3/5/2/1,
24/3/4/5/2/1
}
{
    \newcommand{\sethgt}[1]
    {
        \if##1\a\def\hgt{1}\fi
        \if##1\b\def\hgt{2}\fi
        \if##1\c\def\hgt{3}\fi
        \if##1\d\def\hgt{4}\fi
        \if##1\e\def\hgt{5}\fi
    }
        \foreach\y in {1,...,\n}
        {
            \sethgt{\y} %sets \hgt
            \setclr{\y} %sets \clr
            \fill[\clr] (\i * \stepright, \hgt * \stepup) -- (\i * \stepright + \stepright, \hgt * \stepup) -- (\i * \stepright + \stepright, \hgt * \stepup + \stepup) -- (\i * \stepright, \hgt * \stepup + \stepup) -- cycle;
        }
}

%grid
    \foreach\y in {0,...,\n}
    {
        \draw (\stepright, \y * \stepup + \stepup) -- (\N * \stepright + \stepright, \y * \stepup + \stepup);
    }
    \foreach\x in {0,...,\N}
    {
        \draw (\x * \stepright + \stepright, \stepup) -- (\x * \stepright + \stepright, \n * \stepup + \stepup);
    }   
    
\end{tikzpicture}

$\sigma(L)$
\end{minipage}%
\hfill
\begin{minipage}[t]{.19\textwidth}
\centering
% Almost balanced Gray Code \sigma^{2}(L)
\begin{tikzpicture}[rotate=270]
    \newcommand{\n}{5}                          %size of the permutations
    \newcommand{\width}{10}                       
    \newcommand{\height}{2}                       

    \pgfmathsetmacro\N{24}
    \pgfmathsetmacro\stepright{\width / \N}
    \pgfmathsetmacro\stepup{\height / \n}

%draw colors
\foreach\i/\a/\b/\c/\d/\e in                       %this needs to be expanded if \n is bigger
{1/3/4/5/1/2,
2/3/5/4/1/2,
3/5/3/4/1/2,
4/5/4/3/1/2,
5/4/5/3/1/2,
6/4/3/5/1/2,
7/4/3/1/5/2,
8/1/3/4/5/2,
9/3/1/4/5/2,
10/4/1/3/5/2,
11/1/4/3/5/2,
12/3/4/1/5/2,
13/3/5/1/4/2,
14/1/5/3/4/2,
15/1/3/5/4/2,
16/5/3/1/4/2,
17/5/1/3/4/2,
18/3/1/5/4/2,
19/4/1/5/3/2,
20/1/4/5/3/2,
21/1/5/4/3/2,
22/5/1/4/3/2,
23/5/4/1/3/2,
24/4/5/1/3/2
}
{
    \newcommand{\sethgt}[1]
    {
        \if##1\a\def\hgt{1}\fi
        \if##1\b\def\hgt{2}\fi
        \if##1\c\def\hgt{3}\fi
        \if##1\d\def\hgt{4}\fi
        \if##1\e\def\hgt{5}\fi
    }
        \foreach\y in {1,...,\n}
        {
            \sethgt{\y} %sets \hgt
            \setclr{\y} %sets \clr
            \fill[\clr] (\i * \stepright, \hgt * \stepup) -- (\i * \stepright + \stepright, \hgt * \stepup) -- (\i * \stepright + \stepright, \hgt * \stepup + \stepup) -- (\i * \stepright, \hgt * \stepup + \stepup) -- cycle;
        }
}

%grid
    \foreach\y in {0,...,\n}
    {
        \draw (\stepright, \y * \stepup + \stepup) -- (\N * \stepright + \stepright, \y * \stepup + \stepup);
    }
    \foreach\x in {0,...,\N}
    {
        \draw (\x * \stepright + \stepright, \stepup) -- (\x * \stepright + \stepright, \n * \stepup + \stepup);
    }   
    
\end{tikzpicture}

$\sigma^{2}(L)$
\end{minipage}%
\hfill
\begin{minipage}[t]{.19\textwidth}
\centering
% Almost balanced Gray Code \sigma^{3}(L)
\begin{tikzpicture}[rotate=270]
    \newcommand{\n}{5}                          %size of the permutations
    \newcommand{\width}{10}                       
    \newcommand{\height}{2}                       

    \pgfmathsetmacro\N{24}
    \pgfmathsetmacro\stepright{\width / \N}
    \pgfmathsetmacro\stepup{\height / \n}

%draw colors
\foreach\i/\a/\b/\c/\d/\e in                       %this needs to be expanded if \n is bigger
{1/4/5/1/2/3,
2/4/1/5/2/3,
3/1/4/5/2/3,
4/1/5/4/2/3,
5/5/1/4/2/3,
6/5/4/1/2/3,
7/5/4/2/1/3,
8/2/4/5/1/3,
9/4/2/5/1/3,
10/5/2/4/1/3,
11/2/5/4/1/3,
12/4/5/2/1/3,
13/4/1/2/5/3,
14/2/1/4/5/3,
15/2/4/1/5/3,
16/1/4/2/5/3,
17/1/2/4/5/3,
18/4/2/1/5/3,
19/5/2/1/4/3,
20/2/5/1/4/3,
21/2/1/5/4/3,
22/1/2/5/4/3,
23/1/5/2/4/3,
24/5/1/2/4/3
}
{
    \newcommand{\sethgt}[1]
    {
        \if##1\a\def\hgt{1}\fi
        \if##1\b\def\hgt{2}\fi
        \if##1\c\def\hgt{3}\fi
        \if##1\d\def\hgt{4}\fi
        \if##1\e\def\hgt{5}\fi
    }
        \foreach\y in {1,...,\n}
        {
            \sethgt{\y} %sets \hgt
            \setclr{\y} %sets \clr
            \fill[\clr] (\i * \stepright, \hgt * \stepup) -- (\i * \stepright + \stepright, \hgt * \stepup) -- (\i * \stepright + \stepright, \hgt * \stepup + \stepup) -- (\i * \stepright, \hgt * \stepup + \stepup) -- cycle;
        }
}

%grid
    \foreach\y in {0,...,\n}
    {
        \draw (\stepright, \y * \stepup + \stepup) -- (\N * \stepright + \stepright, \y * \stepup + \stepup);
    }
    \foreach\x in {0,...,\N}
    {
        \draw (\x * \stepright + \stepright, \stepup) -- (\x * \stepright + \stepright, \n * \stepup + \stepup);
    }   
    
\end{tikzpicture}

$\sigma^{3}(L)$
\end{minipage}%
\hfill
\begin{minipage}[t]{.19\textwidth}
\centering
% Almost balanced Gray Code \sigma^{4}(L)
\begin{tikzpicture}[rotate=270]
    \newcommand{\n}{5}                          %size of the permutations
    \newcommand{\width}{10}                       
    \newcommand{\height}{2}                       

    \pgfmathsetmacro\N{24}
    \pgfmathsetmacro\stepright{\width / \N}
    \pgfmathsetmacro\stepup{\height / \n}

%draw colors
\foreach\i/\a/\b/\c/\d/\e in                       %this needs to be expanded if \n is bigger
{1/5/1/2/3/4,
2/5/2/1/3/4,
3/2/5/1/3/4,
4/2/1/5/3/4,
5/1/2/5/3/4,
6/1/5/2/3/4,
7/1/5/3/2/4,
8/3/5/1/2/4,
9/5/3/1/2/4,
10/1/3/5/2/4,
11/3/1/5/2/4,
12/5/1/3/2/4,
13/5/2/3/1/4,
14/3/2/5/1/4,
15/3/5/2/1/4,
16/2/5/3/1/4,
17/2/3/5/1/4,
18/5/3/2/1/4,
19/1/3/2/5/4,
20/3/1/2/5/4,
21/3/2/1/5/4,
22/2/3/1/5/4,
23/2/1/3/5/4,
24/1/2/3/5/4
}
{
    \newcommand{\sethgt}[1]
    {
        \if##1\a\def\hgt{1}\fi
        \if##1\b\def\hgt{2}\fi
        \if##1\c\def\hgt{3}\fi
        \if##1\d\def\hgt{4}\fi
        \if##1\e\def\hgt{5}\fi
    }
        \foreach\y in {1,...,\n}
        {
            \sethgt{\y} %sets \hgt
            \setclr{\y} %sets \clr
            \fill[\clr] (\i * \stepright, \hgt * \stepup) -- (\i * \stepright + \stepright, \hgt * \stepup) -- (\i * \stepright + \stepright, \hgt * \stepup + \stepup) -- (\i * \stepright, \hgt * \stepup + \stepup) -- cycle;
        }
}

%grid
    \foreach\y in {0,...,\n}
    {
        \draw (\stepright, \y * \stepup + \stepup) -- (\N * \stepright + \stepright, \y * \stepup + \stepup);
    }
    \foreach\x in {0,...,\N}
    {
        \draw (\x * \stepright + \stepright, \stepup) -- (\x * \stepright + \stepright, \n * \stepup + \stepup);
    }   
    
\end{tikzpicture}

$\sigma^{4}(L)$
\end{minipage}%

\caption{The balanced Gray code for $\Pall_5$ as created in the proof of Proposition \ref{prop:balancedCode} ($1=$~purple, $2=$~blue, $3=$~green, $4=$~yellow and $5=$~red). $L$ is the almost balanced Gray code shown in Figure \ref{fig:almostbalancedL}, while $\sigma := (1\dotsc 5)$ is the cyclic right rotation of numbers. Another view of this example can be seen in \Cref{fig:Pall5}.}
\end{figure}

It is now straightforward to show by induction that $\Pall_n$ admits a balanced Gray code for odd $n$.

\begin{proposition}\label{prop:balancedOdd}
  Let $n\in \mathbb{N}$ be odd. Then there exists a balanced Gray code \[1\dotsc (n-1)n \xrightarrow{G} 1 \dotsc n (n-1)\] for $\Pall_n$ having $(n-1 \, n)$ as the closing transposition.
\end{proposition}
\begin{proof}
  The case $n=1$ is special, since there is no transposition involved. However, we count the empty list of transpositions as being balanced, since each transposition still appears equally often.
  A balanced Gray code $G$ for $\Pall_3$ is given by
  \[ ((1,3), (1,2), (2,3), (1,3), (1,2)). \]
  Starting at $123$, this yields $123 \xrightarrow{G} 132$, leaving $(2,3)$ as the closing transposition.
  Now assume that there is a balanced Gray code \[1\dotsc (n-1)n \xrightarrow{G} 1\dotsc n(n-1)\] for $\Pall_n$ for odd $n\geq 3$. By \Cref{prop:almostbalancedCode} there is an almost balanced Gray code \[12 \dotsc n+1 \xrightarrow{L} 23\dotsc (n+1) 1\] for $\Pall_{n+1}$ and by \Cref{prop:balancedCode}, $L$ can be used to build a balanced Gray code $H$ for~$\Pall_{n+2}$ which is again of the form \[1\dotsc (n+1)(n+2) \xrightarrow{H} 1\dotsc (n+2)(n+1).\] Hence, there is a balanced Gray code for $\Pall_{2k+1}, k \in \mathbb{N}$, of the claimed form.
 \end{proof}

\subsection{A Balanced Gray Code for Permutations of Even $[n]$} \label{sec:balancedEven}
The idea for creating a balanced Gray code for $\Pall_n$, where $n$ is even, is based on combining two different Gray codes for $\Pall_{n-1}$ while following the construction used in \Cref{prop:balancedCode}. The first Gray code we want to use is the one given by \Cref{prop:balancedOdd}. Constructing the second one is our next goal.

If one compares the starting permutation $12 \dotsc n$ and ending permutation $23 \dotsc 1$ of the Gray code~$L$ from \Cref{prop:almostbalancedCode}, it can be seen that the later is obtained by applying the rotation $(1 \dotsc n)$ to the starting permutation, so a rotation to the left. It turns out that we can also obtain an almost balanced Gray code $R$, where the ending permutation is the starting permutation rotated to the right by one step instead, while satisfying the additional constraint that the transpositions appearing in $L$ and $R$ are the same.

\begin{proposition}\label{prop:almostbalancedCodeRight}
  Let $n\in \mathbb{N}$ be odd and let $12\dotsc (n-1)n \xrightarrow{G} 12\dotsc n (n-1)$ be a balanced Gray code for $\Pall_n$ such that $(n-1, n)$ is the transposition closing the cycle. Then there exists an almost balanced Gray code
  \[12\dotsc (n+1) \xrightarrow{R} (n+1)12 \dotsc n\]
  for $\Pall_{n+1}$, such that the multiset of transpositions appearing in $R$ is the same as in $L$, where $L$ is the almost balanced Gray code from \Cref{prop:almostbalancedCode}.
\end{proposition}
\begin{proof}
  By \Cref{lem:balancedCodeH} we obtain a second balanced Gray code $H$ from $G$ such that
  \[ 12\dotsc (n-1)n \xrightarrow{H} 21 \dotsc (n-1)n \text{ with } c_{+}(H) = \{(1, 2)\} \text{ and } c_{-}(H) = \emptyset. \]

  Let $\tau_{k} := (k \dotsc 1)$ and $\sigma := (1 \dotsc n+1)$ act on $S_{n+1}$. Let then $\bar{G} := \sigma(G)$ and $\bar{H} := \sigma(H)$, meaning that a transposition $(i, j)$ in $G$ (or $H$) becomes $(i+1, j+1)$ in $\bar{G}$ (or $\bar{H}$).
  The diagram shown in \Cref{fig:almostbalancedCodeRight} presents a Gray code for $\Pall_{n+1}$ assembled in a sense reversed from \Cref{prop:almostbalancedCode} such that the end position is a right rotation of the start position. We denote this code by $R$.

  \begin{figure}[h]
    \centering
    \[
      \begin{tikzcd}[row sep=large]
        123\dotsc n(n+1) \arrow[rr, "\bar{G}"]         &        & 123\dotsc (n+1)n \arrow[lld, "{(1, 2)}"']  \\
        213\dotsc(n+1)n \arrow[rr, "\overleftarrow{\tau_{2}(\bar{G})}"]     &        & 213\dotsc n(n+1) \arrow[lld, "{(2, 3)}"']  \\
        312\dotsc n(n+1) \arrow[rr, "\tau_3(\bar{H})"]         &        & 321\dotsc n(n+1) \arrow[lld, "{(3, 4)}"']  \\
        421\dotsc n(n+1) \arrow[rrd, phantom, "\vdots" description] & & \\
        & &  n21\dotsc (n+1) \arrow[lld, "{(n, n+1)}"'] \\
        (n+1)21\dotsc n \arrow[rr, "\overleftarrow{\tau_{n+1}(\bar{H})}"] &        & (n+1)12\dotsc n
      \end{tikzcd}
    \]
    \caption{A second almost balanced Gray code for $\Pall_{n}$ with $n$ even.}
    \label{fig:almostbalancedCodeRight}
  \end{figure}

  Note that $1$ is not contained in any transposition found in $\bar{H}$ or $\bar{G}$. As a result of \Cref{lem:cycleapplication}, the codes $\tau_l(\bar{H})$ and $\tau_l(\bar{G})$ do not change $\tau_l(1) = l$. Hence, permutations of $S_{n+1}$ starting with odd $k$ for $3 \leq k \leq n$ are listed by 
  \[\tau_k(123\dotsc (n+1)) \xrightarrow{\tau_{k}(\bar{H})} \tau_k(132\dotsc (n+1))\] 
  while the ones starting with even $l$ for $4 \leq l \leq n+1$ are given by 
  \[\tau_l(132\dotsc (n+1)) \xrightarrow{\overleftarrow{\tau_{l}(\bar{H})}} \tau_l(123\dotsc (n+1))\]
  The permutations starting with $1$ or $2$ are listed by 
  \[1 \dotsc n(n+1) \xrightarrow{\bar{G}} 1 \dotsc (n+1)n \]
  and 
  \[\tau_2(1 \dotsc (n+1)n) \xrightarrow{\overleftarrow{\tau_{2}(\bar{G})}} \tau_2(1 \dotsc n(n+1))\] 
  respectively. For $k\in\{1, \dotsc, n\}$ the code listing the permutations starting with $k$ is always connected by the transposition~$(k, k+1)$ to the code listing the permutations starting with $k+1$.

  Proving $R$ to be almost balanced is similar to \Cref{prop:almostbalancedCode}. Fix ${k \in \{3,\dotsc,n+1\}}$. Because $H$ is balanced, each transposition on $S_{n+1}$ not containing $1$ appears $2(n-2)!$ times in $H$, except for~$(2, 3)$ which appears $2(n-2)! - 1$ times. Since $\tau_{k}$ is bijective such that $\tau_{k}((2, 3)) = (1, 2)$ and $\tau_k(H)$ leaves~$k$ unchanged,
  it follows that all transpositions on~$S_{n+1}$ not containing $k$ appear $2(n-2)!$ times in $\tau_{k}(H)$ except for $(1, 2)$.

  Focusing on $G$, we obtain that all transpositions on $S_{n+1}$ not containing $1$ appear ${2(n-2)!}$ times in $G$ except for $(n, n+1)$. Lastly, $\tau_{2}(G)$ contains all transpositions on $S_{n+1}$ not containing $2$ exactly $2(n-2)!$ times except for $(n, n+1)$.

  Hence, for $(i, j)$ a wide transposition on $S_{n+1}$, we have shown that $(i, j)$ appears
  \begin{itemize}
    \item $2(n-2)!$ times in $\tau_{k}(H)$ for $k\in \{3,\dotsc,n+1\}\setminus \{i,j\}$ and
    \item $2(n-2)!$ times in $\tau_{k}(G)$ for $k\in \{1,2\}\setminus \{i,j\}$.
  \end{itemize}
  Consequently, $(i, j)$ is contained $(n-1) \cdot 2(n-2)! = 2(n-1)!$ times in $R$ making it almost balanced.

  Let $L$ be the almost balanced Gray code from \Cref{prop:almostbalancedCode} and define ${\omega_{k} := (k \dotsc n+1)}$ acting on $S_{n+1}$. To show that the multiset of transpositions in $R$ is the same as in $L$ observe that for $k\in \{3,\dotsc,n+1\}$ the mulitset of transpositions in $\tau_{k}(H)$ is the same as in $\omega_{k}(G)$,
  since all transpositions on $S_{n+1}$ not containing $k$ appear $2(n-2)!$ times within each except for $(1, 2)$ which appears $2(n-2)! - 1$ times in both.
  The same holds for~$k\in \{1,2\}$ just that $(n, n+1)$ is the transposition appearing only $2(n-2)! - 1$ times in both. Finally, the set of transpositions connecting these partial Gray codes is the same in $L$ and $R$, namely $\{(i, i+1) \mid i\in [n] \}$. As a result, the claim follows.
\end{proof}

For odd $n\in \mathbb{N}$ we are now in the position to construct the second special Gray code for~$\Pall_n$. Note that this code is neither balanced nor cyclic.

\begin{proposition}\label{prop:secondGrayCode}
  Let $n > 1$ be odd and denote by $G$ the balanced Gray code obtained for $\Pall_n$ from \Cref{prop:balancedOdd}. There exists a Gray code $H$ for $\Pall_n$ such that
  \[ 23\dotsc (n-1)1n \xrightarrow{H} n12 \dotsc (n-1)\]
  with the additional property that the multiset of transpositions in $G$ and in $H$ is the same.
\end{proposition}
\begin{proof}
  Denote by $L$ the almost balanced Gray code for $\Pall_{n-1}$ obtained from \Cref{prop:almostbalancedCode} and by $R$ the almost balanced Gray code constructed in \Cref{prop:almostbalancedCodeRight}. Further, let $\sigma := (1 2 \dotsc n)$ act on $S_{n}$. We present the diagram for the Gray code $H$ in \Cref{fig:secondGrayCode}.

  \begin{figure}[h]
    \centering
    \[
      \begin{tikzcd}[row sep=large]
        23\dotsc (n-1)1n \arrow[rr, "\overleftarrow{L}"]                      &        & 12 \dotsc n \arrow[lld, "{(1, n)}"']             \\
        n23\dotsc (n-1)1 \arrow[rr, "\overleftarrow{\sigma(R)}"]                      &        & 23 \dotsc (n-1)n1 \arrow[lld, "{(1, 2)}"']       \\
        13\dotsc (n-1)n2 \arrow[rr, "\overleftarrow{\sigma^{2}(R)}"]              &        & 34 \dotsc n12 \arrow[lld, "{(2, 3)}"']           \\
        24\dotsc n13 \arrow[rrd, phantom, "\vdots" description] & & \\
        & & (n-1)n1\dotsc (n-2) \arrow[lld, "{(n-2, n-1)}"']\\
        (n-2)n1\dotsc (n-1) \arrow[rr, "\overleftarrow{\sigma^{n-1}(R)}"] &        & n12 \dotsc (n-1)
      \end{tikzcd}
    \]
    \caption{Gray code for $\Pall_{n}$ with $n$ odd.}
    \label{fig:secondGrayCode}
  \end{figure}

  The permutations ending with $n$ are listed by the partial Gray code 
  \[23\dotsc (n-1)1n \xrightarrow{\overleftarrow{L}} 12 \dotsc n\] 
  and for $k \in \{1,\dots, n-1\}$ the permutations ending with $k$ are listed by 
  \[ \sigma^k((n-1)1\dotsc (n-2)n) \xrightarrow{\overleftarrow{\sigma^{k}(R)}} \sigma^k(12\dotsc (n-1)n).\] 
  Consequently, the transposition connecting $\overleftarrow{\sigma^{k}(R)}$ and $\overleftarrow{\sigma^{k+1}(R)}$ is $(k, k+1)$. The code $\overleftarrow{L}$ is connected to $\overleftarrow{\sigma(R)}$ by $(1, n)$.

  For the rest of this proof we treat the codes $L, R, H$ and $G$ as multisets of transpositions by forgetting the order. We thus have
  \[ H = \bigcupdot_{k\in [n-1]}\sigma^{k}(R) \cupdot L \cupdot \{\sigma^{k}((1, n)) \mid k \in \{0,\dotsc,n-2\}\}. \]
  In order to prove $H = G$, we use the fact that $L$ and $R$ contain the same transpositions by \Cref{prop:almostbalancedCodeRight}. Thus, since $\sigma$ is bijective,
  $\sigma^{k}(L) = \sigma^{k}(R)$.
  This already yields
  \[H = \bigcupdot_{k\in \{0,\dotsc,n-1\}}\sigma^{k}(L) \cupdot \{\sigma^{k}((1, n)) \mid k \in \{0,\dotsc,n-2\}\} = G.\]
\end{proof}

Finally, let us combine the obtained Gray codes for $\Pall_{n-1}$ with even $n\in\mathbb{N} \setminus \{0\}$ into a balanced Gray code for $\Pall_n$.

\begin{proposition}\label{prop:balancedEven}
  Let $n\in \mathbb{N} \setminus \{0\}$ be even. Then there exists a balanced Gray code $G$
  \[ 12 \dotsc n \xrightarrow{G} n2 \dotsc (n-1)1 \]
  for $\Pall_n$ such that $(1, n)$ is the transposition closing the cycle.
\end{proposition}
\begin{proof}
  The case $n=2$ is trivial. Hence, we assume $4 \leq n$.
  Let $H_{1}$ be the balanced Gray code obtained from \Cref{prop:balancedOdd} for $\Pall_{n-1}$ and denote by $H_{2}$ the Gray code for~$\Pall_{n-1}$ constructed in \Cref{prop:secondGrayCode}. Further let $\sigma := (n \dotsc 2 1)$ act on $S_{n-1}$. The diagram in \Cref{fig:balancedEven} presents the cyclic Gray code $G$ which we want to prove to be balanced for $\Pall_n$.

  \begin{figure}[h]
    \centering
    \[
      \begin{tikzcd}[row sep=large, column sep=large]
        1\dotsc (n-2)(n-1)n \arrow[rr, "H_{1}"]                                                                  &        & 1\dotsc (n-1)(n-2)n \arrow[lld, "{(n-1, n)}"']            \\
        1\dotsc n(n-2)(n-1) \arrow[rr, "\sigma(H_{2})"]  &        & (n-2)n1 \dotsc (n-1) \arrow[lld, "{(n-2, n-1)}"']     \\
        (n-1)n1\dotsc (n-2)  \arrow[rrd, phantom, "\vdots" description]   &        & {} \\
        &  & 3 \dotsc (n-1)1n2 \arrow[lld, "(1 2)"']                                     \\
        3 \dotsc (n-1)2n1 \arrow[rr, "\sigma^{n-1}(H_{2})"]                                                   &        & n23\dotsc (n-1)1
      \end{tikzcd}
    \]
    \caption{Balanced Gray code for $\Pall_{n}$ with $n$ even.}
    \label{fig:balancedEven}
  \end{figure}

  The permutations of $S_{n}$ ending with even $k \in [n]$ are listed by 
  \[ C_k := \sigma^{k}(1 \dotsc (n-3)(n-2)(n-1)n) \xrightarrow{\sigma^{k}(H_{1})} \sigma^{k}(1\dotsc (n-3)(n-1)(n-2)n)\]
  and the ones ending with odd $l\in [n]$ by 
  \[ D_l := \sigma^{l}(2\dotsc (n-2)1(n-1)n) \xrightarrow{\sigma^{l}(H_{2})} \sigma^{l}((n-1)1\dotsc (n-2)n).\] 
  To see that $G$ is well defined, observe that the ending permutation of $C_0$ differs by the transposition $(n-1, n)$ from the starting permutation of $D_{1}$. As a result, the ending permutation of $C_k$ differs by the transposition $\sigma^{k}((n-1, n))$ from the starting permutation of $D_{k+1}$. 
  On the other hand, the ending permutation of $D_0$ also differs by the transposition~$(n-1, n)$ from the starting permutation of $C_1$. Consequently, applying~$\sigma^{k}((n-1, n))$ to the ending permutation of $D_k$ yields the starting permutation of $C_{k+1}$. Thus, the Gray code is well-defined. Next, we prove $G$ to be balanced.

  In the following we treat the codes $H_1, H_2$ and $G$ as multisets of transpositions by forgetting the order. Define
\[Z := \{\sigma^{k}((1, 2)) \mid k \in \{0,\dotsc,n-2\}\} = \{(i, i+1) \mid i \in [n-1]\}. \]
  Note that by \Cref{prop:secondGrayCode} we have $H_{1} = H_{2}$ and, therefore,
  \[G = \bigcupdot_{k \in \{0,\dotsc,n-1\}}\sigma^{k}(H_{1}) \cupdot Z. \]
  No transposition in $H_{1}$ contains $n$ and, hence, by \Cref{lem:cycleapplication} no transposition in $\sigma^{k}(H_{1})$ contains $n-k$. Since $\sigma$ is bijective and $H_{1}$ is balanced such that ${c_{+}(H_{1}) = \{(1, 2)\}}$, it follows that each transposition on $S_{n}$ not containing $n-k$ appears $2(n-2)!$ times in $\sigma^{k}(H_{1})$ except for $\sigma^{k}((1, 2))$ which appears ${2(n-2)! - 1}$ times.
  It then holds that a transposition on $S_{n}$ not contained in $Z \cupdot \{(1, n)\}$ appears $(n-1) \cdot 2(n-2)! = 2(n-1)!$ times in $G$. A transposition from $Z$, on the other hand, appears \[(n-1) \cdot (2(n-2)! - 1) + (n-2) = 2(n-1)! - 1\] times in $\bigcupdot_{k \in \{0,\dotsc,n-1\}} \sigma^{k}(H_{1})$ and once in $Z$. As a result, it is also balanced in $G$. The only remaining transposition is $(1, n)$. This one also appears $2(n-1)! - 1$ times in $G$, like the transpositions from $Z$. Additionally, it is the transposition closing the cycle. Hence, $G$ is a balanced Gray code for $\Pall_n$.
\end{proof}

We want to demonstrate the construction for $\Pall_4$.
\begin{example}[Balanced Gray code for $\Pall_4$] \label{ex:balancedEven4}
  Denote by \[H_{1} := 123 \xrightarrow{((1,3), (1,2), (2,3), (1,3), (1,2))} 132 \] the balanced Gray code from \Cref{prop:balancedOdd} and by
  \[ H_{2} := 213 \xrightarrow{((1,2), (1,3), (2,3), (1,2), (1,3))} 312\] the code from \Cref{prop:secondGrayCode}. Further, let $\sigma := (4 \, 3 \, 2 \, 1)$ act on $S_{4}$. The blue edges in the following diagram show the resulting balanced Gray code $G$ for $\Pall_4$, following the construction described in \Cref{prop:balancedEven}, while \Cref{fig:balancedEven4} depicts the resulting code visually.
  \[
    \begin{tikzcd}[row sep=huge, column sep=huge]
      1234 \arrow[r, "H_1"', blue]           & 1324 \arrow[rd, "{(3, 4)}", blue] & 2134 \arrow[r, "H_2"]           & 3124                             \\
      4123 \arrow[r, "\sigma(H_1)"'] & 4213                          & 1423 \arrow[r, "\sigma(H_2)", blue] & 2413 \arrow[llld, "{(2, 3)}"', blue] \\
      3412 \arrow[r, "\sigma^{2}(H_1)"', blue]   & 3142 \arrow[rd, "{(1, 2)}", blue] & 4312 \arrow[r, "\sigma^{2}(H_2)"]   & 1342                             \\
      2341 \arrow[r, "\sigma^{3}(H_1)"']           & 2431                          & 3241 \arrow[r, "\sigma^{3}(H_2)", blue]          & 4231
    \end{tikzcd}
  \]

    %
% Creation of balanced Gray Code for S_4
%

\begin{figure}[h!]
\begin{minipage}[t]{.24\textwidth}
\centering
% H_1 4
\begin{tikzpicture}[rotate=270]
    \newcommand{\n}{4}                          %size of the permutations
    \newcommand{\width}{3}                       
    \newcommand{\height}{2}  

    \pgfmathsetmacro\N{6}
    \pgfmathsetmacro\stepright{\width / \N}
    \pgfmathsetmacro\stepup{\height / \n}

%draw colors
\foreach\i/\a/\b/\c/\d in                       %this needs to be expanded if \n is bigger
{1/1/2/3/4,
2/3/2/1/4,
3/3/1/2/4,
4/2/1/3/4,
5/2/3/1/4,
6/1/3/2/4
}
{
    \newcommand{\sethgt}[1]
    {
        \if##1\a\def\hgt{1}\fi
        \if##1\b\def\hgt{2}\fi
        \if##1\c\def\hgt{3}\fi
        \if##1\d\def\hgt{4}\fi
    }
        \foreach\y in {1,...,\n}
        {
            \sethgt{\y} %sets \hgt
            \setclr{\y} %sets \clr
            \fill[\clr] (\i * \stepright, \hgt * \stepup) -- (\i * \stepright + \stepright, \hgt * \stepup) -- (\i * \stepright + \stepright, \hgt * \stepup + \stepup) -- (\i * \stepright, \hgt * \stepup + \stepup) -- cycle;
        }
}

%grid
    \foreach\y in {0,...,\n}
    {
        \draw (\stepright, \y * \stepup + \stepup) -- (\N * \stepright + \stepright, \y * \stepup + \stepup);
    }
    \foreach\x in {0,...,\N}
    {
        \draw (\x * \stepright + \stepright, \stepup) -- (\x * \stepright + \stepright, \n * \stepup + \stepup);
    }   
    
\end{tikzpicture}

$H_1$
\end{minipage}%
\hfill
\begin{minipage}[t]{.24\textwidth}
\centering
% \sigma(H_2 4)
\begin{tikzpicture}[rotate=270]
    \newcommand{\n}{4}                          %size of the permutations
    \newcommand{\width}{3}                       
    \newcommand{\height}{2}  

    \pgfmathsetmacro\N{6}
    \pgfmathsetmacro\stepright{\width / \N}
    \pgfmathsetmacro\stepup{\height / \n}

%draw colors
\foreach\i/\a/\b/\c/\d in                       %this needs to be expanded if \n is bigger
{1/1/4/2/3,
2/4/1/2/3,
3/2/1/4/3,
4/1/2/4/3,
5/4/2/1/3,
6/2/4/1/3
}
{
    \newcommand{\sethgt}[1]
    {
        \if##1\a\def\hgt{1}\fi
        \if##1\b\def\hgt{2}\fi
        \if##1\c\def\hgt{3}\fi
        \if##1\d\def\hgt{4}\fi
    }
        \foreach\y in {1,...,\n}
        {
            \sethgt{\y} %sets \hgt
            \setclr{\y} %sets \clr
            \fill[\clr] (\i * \stepright, \hgt * \stepup) -- (\i * \stepright + \stepright, \hgt * \stepup) -- (\i * \stepright + \stepright, \hgt * \stepup + \stepup) -- (\i * \stepright, \hgt * \stepup + \stepup) -- cycle;
        }
}

%grid
    \foreach\y in {0,...,\n}
    {
        \draw (\stepright, \y * \stepup + \stepup) -- (\N * \stepright + \stepright, \y * \stepup + \stepup);
    }
    \foreach\x in {0,...,\N}
    {
        \draw (\x * \stepright + \stepright, \stepup) -- (\x * \stepright + \stepright, \n * \stepup + \stepup);
    }   
    
\end{tikzpicture}

$\sigma(H_2)$
\end{minipage}%
\hfill
\begin{minipage}[t]{.24\textwidth}
\centering
% \sigma^{2}(H_1 4)
\begin{tikzpicture}[rotate=270]
    \newcommand{\n}{4}                          %size of the permutations
    \newcommand{\width}{3}                       
    \newcommand{\height}{2}                       

    \pgfmathsetmacro\N{6}
    \pgfmathsetmacro\stepright{\width / \N}
    \pgfmathsetmacro\stepup{\height / \n}

%draw colors
\foreach\i/\a/\b/\c/\d in                       %this needs to be expanded if \n is bigger
{1/3/4/1/2,
2/1/4/3/2,
3/1/3/4/2,
4/4/3/1/2,
5/4/1/3/2,
6/3/1/4/2
}
{
    \newcommand{\sethgt}[1]
    {
        \if##1\a\def\hgt{1}\fi
        \if##1\b\def\hgt{2}\fi
        \if##1\c\def\hgt{3}\fi
        \if##1\d\def\hgt{4}\fi
    }
        \foreach\y in {1,...,\n}
        {
            \sethgt{\y} %sets \hgt
            \setclr{\y} %sets \clr
            \fill[\clr] (\i * \stepright, \hgt * \stepup) -- (\i * \stepright + \stepright, \hgt * \stepup) -- (\i * \stepright + \stepright, \hgt * \stepup + \stepup) -- (\i * \stepright, \hgt * \stepup + \stepup) -- cycle;
        }
}

%grid
    \foreach\y in {0,...,\n}
    {
        \draw (\stepright, \y * \stepup + \stepup) -- (\N * \stepright + \stepright, \y * \stepup + \stepup);
    }
    \foreach\x in {0,...,\N}
    {
        \draw (\x * \stepright + \stepright, \stepup) -- (\x * \stepright + \stepright, \n * \stepup + \stepup);
    }   
    
\end{tikzpicture}

$\sigma^{2}(H_1)$
\end{minipage}%
\hfill
\begin{minipage}[t]{.24\textwidth}
\centering
% \sigma^{3}(H_2 4)
\begin{tikzpicture}[rotate=270]
    \newcommand{\n}{4}                          %size of the permutations
    \newcommand{\width}{3}                       
    \newcommand{\height}{2}                       

    \pgfmathsetmacro\N{6}
    \pgfmathsetmacro\stepright{\width / \N}
    \pgfmathsetmacro\stepup{\height / \n}

%draw colors
\foreach\i/\a/\b/\c/\d in                       %this needs to be expanded if \n is bigger
{1/3/2/4/1,
2/2/3/4/1,
3/4/3/2/1,
4/3/4/2/1,
5/2/4/3/1,
6/4/2/3/1
}
{
    \newcommand{\sethgt}[1]
    {
        \if##1\a\def\hgt{1}\fi
        \if##1\b\def\hgt{2}\fi
        \if##1\c\def\hgt{3}\fi
        \if##1\d\def\hgt{4}\fi
    }
        \foreach\y in {1,...,\n}
        {
            \sethgt{\y} %sets \hgt
            \setclr{\y} %sets \clr
            \fill[\clr] (\i * \stepright, \hgt * \stepup) -- (\i * \stepright + \stepright, \hgt * \stepup) -- (\i * \stepright + \stepright, \hgt * \stepup + \stepup) -- (\i * \stepright, \hgt * \stepup + \stepup) -- cycle;
        }
}

%grid
    \foreach\y in {0,...,\n}
    {
        \draw (\stepright, \y * \stepup + \stepup) -- (\N * \stepright + \stepright, \y * \stepup + \stepup);
    }
    \foreach\x in {0,...,\N}
    {
        \draw (\x * \stepright + \stepright, \stepup) -- (\x * \stepright + \stepright, \n * \stepup + \stepup);
    }   
    
\end{tikzpicture}

$\sigma^{3}(H_2)$
\end{minipage}%

\caption{Visualization of balanced Gray code for $\Pall_4$ as obtained in Example \ref{ex:balancedEven4} where $1=$~purple, $2=$~blue, $3=$~green and $4=$~yellow and $\sigma := (4\, 3\, 2\, 1)$.}
\label{fig:balancedEven4}
\end{figure}

\end{example}

We are finally in the position to prove \Cref{thm:perm_rainbow}.

\begin{proof}[Proof of \Cref{thm:perm_rainbow}] \label{proof:theorem2}
  Let $n \in \mathbb{N}^{+}$, fix $m \in [n]$ and define $M := 2(m-2)!$. By using either \Cref{prop:balancedOdd} or \Cref{prop:balancedEven}, depending on the parity of $m$, there exists a balanced Gray code for $\Pall_m$. We will show that if $G$ is a $M$-rainbow cycle in~$\Pall_k$ for $m \leq k < n$, then we can construct a $M$-rainbow cycle in $\Pall_{k+1}$. The result then follows by induction.

  So fix $k$ with $m \leq k < n$ and define $N := M \cdot \binom{k}{2}$. Let $G := (t_{i})_{i\in [N - 1]}$ be a $M$-rainbow cycle in~$\Pall_k$.
  For a permutation ${\sigma = \sigma_{1}\dotsc\sigma_{k} \in S_{k}}$ denote the partial Gray code
  \[((\sigma_{k}, k+1), (\sigma_{k-1}, k+1), \dotsc, (\sigma_{1}, k+1)) \]
  pushing $k+1$ from right to left on $\sigma$ by $L(\sigma)$. Note that each transposition on $S_{k+1}$ containing $k+1$ appears once in $L(\sigma)$.
  We can then construct a $M$-rainbow cycle in $\Pall_{k+1}$ starting on an arbitrary $\pi_1 \in S_{k+1}$, defining $ \pi_i \xrightarrow{t_i} \pi_{i+1} $, as follows:
  
  \[
     \begin{tikzcd}
       \pi_1(k+1) \arrow[r, "L(\pi_1)"] & (k+1)\pi_1 \arrow[r, "t_1"] & (k+1)\pi_2 \arrow[r, "\overleftarrow{L(\pi_2)}"] & \dotsc \arrow[r, "\overleftarrow{L(\pi_{M})}"] & \pi_{M}(k+1) \arrow[d, "t_{M}"]         \\
       \pi_N(k+1)                       & \dots \arrow[l, "t_{N-1}"'] & \pi_{M+2}(k+1) \arrow[l, "t_{M+2}"']             &                                               & \pi_{M + 1}(k+1) \arrow[ll, "t_{M+1}"']
     \end{tikzcd}
 \]
 The trick is to keep pushing $k+1$ from right to left on $\pi_{i}$ and from left to right on $\pi_{i+1}$, for $1 \leq i < M$ odd, a total number of $M$ times. This implies that each transposition on~$S_{k+1}$ containing $k+1$ appears~$M$ times. By assumption all transpositions on $S_{k+1}$ not containing $k+1$ already appear $M$ times in~$G$.
\end{proof}

\subsection{Balancing Cyclically Adjacent Transpositions} \label{sec:cyclicallyBalanced}
In this section we prove \Cref{thm:cyclicallyballanced}.

\cyclicallyBalanced*

In this section, a transposition will not be applied to the entries but to the indices of a permutation. Recall that the cyclically adjacent transpositions on $S_{n}$ are $(i, i+1)$ for~$i \in [n-1]$ and $(1, n)$. A balanced Gray code for $\Pcadj_n$ is, hence, a cyclic listing of $S_{n}$ such that each of those transpositions appears exactly $(n-1)!$ times.
We will treat only the case where $n$ is even, since for odd $n$ the construction is given in Theorem 5.7 in \cite{Gregor2023}.

Let $n\in \mathbb{N} \setminus \{0\}$ be even and let $L$ be a Gray code for $\Padj_n$ starting at $1\dots n$ and ending at $23 \dots n1$ given by the fact that the Permutahedron has a Hamilton path between those permutations, if $n$ is even \cite{Tchuente1982}. Starting the code $L$ not at $1 \dotsc n$ but at an arbitrary permutation $\pi_{1}\dots \pi_{n}$, leads to the ending permutation $\pi_{2}\pi_{3}\dots \pi_{n}\pi_{1}$, since transpositions are applied to indices. Furthermore, we can use this code backwards, leading to a Gray code from $\pi_{1}\dots \pi_{n}$ to $\pi_{n}\pi_{1}\dots \pi_{n-1}$ while having the same transpositions as $L$. We denote this code by $R$.

Using $L$ and $R$ we construct two different Gray codes for $\Padj_{n+1}$. These will then be combined to a balanced Gray code for $\Pcadj_{n+2}$. This approach mimics the way how the balanced Gray code was obtained in \Cref{prop:balancedEven}.

\begin{proposition}\label{prop:cyclical_G}
  Let $n\in \mathbb{N} \setminus \{0\}$ be even. Then there is a Gray code for $\Padj_{n+1}$ starting from $12 \dots n+1$ to $21 \dots n+1$. We denote this code by $G$.
\end{proposition}
\begin{proof}
  Let $\sigma := (1 \dotsc n+1)$. \Cref{fig:cyclical_G} shows the desired Gray code. Note that $\sigma(L)$ always keeps the first entry of a permutation constant, while the later $n$ coordinates follow the Gray code $L$.
  \begin{figure}[!ht]
    \centering
      \[
        \begin{tikzcd}[row sep=small]
          12\dots n+1 \arrow[rr, "\sigma(L)"]      &        & 13\dots (n+1)2 \arrow[lld, "{(1, 2)}"'] \\
          31\dots(n+1)2 \arrow[rr, "\sigma(L)"] &        & 34\dots(n+1)21 \arrow[lld, "{(1, 2)}"'] \\
          43\dots(n+1)21                    & \vdots & (n+1)21\dots n \arrow[lld, "{(1, 2)}"'] \\
          2(n+1)1\dots n \arrow[rr, "\sigma(L)"]   &        & 21\dots n+1
        \end{tikzcd}
      \]
  \caption{First Gray code for $\Padj_{n+1}$}
  \label{fig:cyclical_G}
  \end{figure}
  We repeatedly apply $\sigma(L)$ followed by a transposition of the first two coordinates until all permutations of $S_{n+1}$ are listed. 
\end{proof}

Let us move to the second code.

\begin{proposition}\label{prop:cyclical_H}
  Let $n\in \mathbb{N} \setminus \{0\}$ be even. Then there is a Gray code for $\Padj_{n+1}$ starting from $13\dots (n+1)2$ to $(n+1)12\dots n$ having the same transpositions as $G$. We denote this code by $H$.
\end{proposition}
\begin{proof}
  Let $\sigma := (1 \dotsc n+1)$.
  \Cref{fig:cyclical_H} shows the desired Gray code. Note, as in \Cref{prop:cyclical_G}, that~$\sigma(L)$ and $\sigma(R)$ always keep the first entry of a permutation constant, while the later $n$ coordinates follow the Gray code $L$ or $R$ respectively.
  \begin{figure}[h]
    \centering
  \[
    \begin{tikzcd}[row sep=small]
      134\dots (n+1)2 \arrow[rr, "\sigma(R)"]       &        & 12\dots (n+1) \arrow[lld, "{(1, 2)}"']     \\
      21\dots(n+1) \arrow[rr, "\sigma(L)"]          &        & 23\dots(n+1)1 \arrow[lld, "{(1, 2)}"']     \\
      32\dots(n+1)1                          & \vdots & n(n+1)1\dots(n-1) \arrow[lld, "{(1, 2)}"'] \\
      (n+1)n1\dots(n-1) \arrow[rr, "\sigma(L)"] &        & (n+1)1\dots n
    \end{tikzcd}
  \]
  \caption{Second Gray code for $\Padj_{n+1}$}
  \label{fig:cyclical_H}
  \end{figure}
  After initially applying $\sigma(R)$ to $134\dots (n+1)2$ followed by a transposition of the first two entries, we continue by repeatedly applying $\sigma(L)$ followed by a transposition of the first two coordinates until all permutations of $S_{n+1}$ are listed.

  Since the transpositions of $L$ and $R$ are the same, it is clear that the transpositions of $G$ and $H$ are the same.
\end{proof}

We are now in the position to prove \Cref{thm:cyclicallyballanced}.

\begin{proof}[Proof of \Cref{thm:cyclicallyballanced}]
\label{proof:theorem3}
  Let $n\geq 1$ be odd and let $\sigma := (1 \dotsc n+1)$ be the automorphism of~$\Pcadj_{n+1}$ that cyclically shifts all entries one coordinate to the right. Taking $G$ and $H$ from \Cref{prop:cyclical_G} and \Cref{prop:cyclical_H} for $\Padj_n$, we have the following starting and ending permutations for $\sigma^k(G)$ and $\sigma^k(H)$:
  \[ \sigma^k(12\dotsc (n+1)) \xrightarrow{\sigma^k(G)} \sigma^k(21 \dotsc (n+1))\]
  and
  \[ \sigma^k(13\dotsc n2(n+1)) \xrightarrow{\sigma^k(H)} \sigma^k(n1\dotsc (n-1)(n+1)). \]
  \enlargethispage{\baselineskip}
  Furthermore, it then holds that
  \[
    \begin{tikzcd}
    12\dots (n+1) \arrow[rr, "G"]             &  & 21\dots (n+1) \arrow[lld, "{(1, n+1)}" description]      \\
    (n+1)1\dots n2 \arrow[rr, "\sigma(H)"] &  & (n+1)n1\dots (n-1)                                       \\
    \end{tikzcd}
  \]
  and
  
  \[
    \begin{tikzcd}
    13\dots n2(n+1) \arrow[rr, "H"]       &  & n1\dots (n-1)(n+1) \arrow[lld, "{(1, n+1)}" description] \\
    (n+1)1\dots n \arrow[rr, "\sigma(G)"] &  & (n+1)21\dots n.                                      
    \end{tikzcd}
  \]
  
  Observe, that by (an appropriate version of) \Cref{lem:cycleapplication}, $\sigma^k(G)$ and $\sigma^k(H)$ both contain no transposition containing $\sigma^k(1)$, since the coordinate $n+1$ is constant in $G$ and $H$.
  It follows that $\sigma^k(G)$ can be connected to $\sigma^{k+1}(H)$ by the transposition $\sigma^k((1, n+1))$ such that the entry $n+1$ stays constantly at coordinate~$\sigma^k(1)$ in $\sigma^k(G)$ and at coordinate $\sigma^{k+1}(1)$ in $\sigma^{k+1}(H)$.
  On the other hand,~$\sigma^k(H)$ can also be connected to $\sigma^{k+1}(G)$ by the transposition $\sigma^k((1, n+1))$ such that the entry $n+1$ stays constantly at coordinate $\sigma^k(1)$ in $\sigma^k(H)$ and at coordinate $\sigma^{k+1}(1)$ in $\sigma^{k+1}(G)$.
  
  As a result and because $G$ and $H$ are Gray codes for $\Padj_{n}$, the code  
  \[ G, (1,n+1), \sigma(H), (1,2), \sigma^{2}(G), (2,3),\sigma^{3}(H),\dots, \sigma^{n}(H) \]
  yields a cyclic Gray code $C$ for $\Pcadj_{n+1}$, since the entry $n+1$ stays constant within each code and takes each coordinate exactly once.
  
  By construction each transposition in $G$ is adjacent, so it is of the form $(i, i+1)$ for some~$i\in [n-1]$. Hence, it becomes 
  $\sigma^k((i, i+1))$ in $\sigma^k(G)$. Because the transpositions in $\sigma^k(G)$ and
  $\sigma^k(H)$ are equal by \Cref{prop:cyclical_H}, it follows that $C$ consists of the transpositions
  \[ \bigcupdot_{t\in G} \{\sigma^k(t) \mid k\in [n]\} \cupdot \{\sigma^k((1, n+1)) \mid k\in [n]\}\]
  where the later set refers to the transpositions connecting the partial codes including the one closing the cycle.
  Note that 
  \[ \bigcupdot_{t\in G} \{\sigma^k(t) \mid k\in [n]\} = \bigcupdot_{i\in [n! - 1]} \{\sigma^k((1,n+1)) \mid k\in [n]\}.\]
  Consequently, each cyclically adjacent transposition on $S_{n+1}$ appears exactly $n!$ times in~$C$.
\end{proof}

We end this section by giving the balanced Gray code for $\Pcadj_4$ in \Cref{fig:balancedAdjacentExample}.
%
% Creation of balanced Gray Code for S_4
%

\begin{figure}[h!]
\centering
\begin{minipage}[t]{.24\textwidth}
\centering
% G
\begin{tikzpicture}[rotate=270]
    \newcommand{\n}{4}                          %size of the permutations
    \newcommand{\width}{3}                       
    \newcommand{\height}{2}  

    \pgfmathsetmacro\N{6}
    \pgfmathsetmacro\stepright{\width / \N}
    \pgfmathsetmacro\stepup{\height / \n}

%draw colors
\foreach\i/\a/\b/\c/\d in                       %this needs to be expanded if \n is bigger
{1/1/2/3/4,
2/1/3/2/4,
3/3/1/2/4,
4/3/2/1/4,
5/2/3/1/4,
6/2/1/3/4
}
{
    \newcommand{\sethgt}[1]
    {
        \if##1\a\def\hgt{1}\fi
        \if##1\b\def\hgt{2}\fi
        \if##1\c\def\hgt{3}\fi
        \if##1\d\def\hgt{4}\fi
    }
        \foreach\y in {1,...,\n}
        {
            \sethgt{\y} %sets \hgt
            \setclr{\y} %sets \clr
            \fill[\clr] (\i * \stepright, \hgt * \stepup) -- (\i * \stepright + \stepright, \hgt * \stepup) -- (\i * \stepright + \stepright, \hgt * \stepup + \stepup) -- (\i * \stepright, \hgt * \stepup + \stepup) -- cycle;
        }
}

%grid
    \foreach\y in {0,...,\n}
    {
        \draw (\stepright, \y * \stepup + \stepup) -- (\N * \stepright + \stepright, \y * \stepup + \stepup);
    }
    \foreach\x in {0,...,\N}
    {
        \draw (\x * \stepright + \stepright, \stepup) -- (\x * \stepright + \stepright, \n * \stepup + \stepup);
    }   
    
\end{tikzpicture}

$G$
\end{minipage}%
\hfill
\begin{minipage}[t]{.24\textwidth}
\centering
% \sigma(H)
\begin{tikzpicture}[rotate=270]
    \newcommand{\n}{4}                          %size of the permutations
    \newcommand{\width}{3}                       
    \newcommand{\height}{2}  

    \pgfmathsetmacro\N{6}
    \pgfmathsetmacro\stepright{\width / \N}
    \pgfmathsetmacro\stepup{\height / \n}

%draw colors
\foreach\i/\a/\b/\c/\d in                       %this needs to be expanded if \n is bigger
{1/4/1/3/2,
2/4/1/2/3,
3/4/2/1/3,
4/4/2/3/1,
5/4/3/2/1,
6/4/3/1/2
}
{
    \newcommand{\sethgt}[1]
    {
        \if##1\a\def\hgt{1}\fi
        \if##1\b\def\hgt{2}\fi
        \if##1\c\def\hgt{3}\fi
        \if##1\d\def\hgt{4}\fi
    }
        \foreach\y in {1,...,\n}
        {
            \sethgt{\y} %sets \hgt
            \setclr{\y} %sets \clr
            \fill[\clr] (\i * \stepright, \hgt * \stepup) -- (\i * \stepright + \stepright, \hgt * \stepup) -- (\i * \stepright + \stepright, \hgt * \stepup + \stepup) -- (\i * \stepright, \hgt * \stepup + \stepup) -- cycle;
        }
}

%grid
    \foreach\y in {0,...,\n}
    {
        \draw (\stepright, \y * \stepup + \stepup) -- (\N * \stepright + \stepright, \y * \stepup + \stepup);
    }
    \foreach\x in {0,...,\N}
    {
        \draw (\x * \stepright + \stepright, \stepup) -- (\x * \stepright + \stepright, \n * \stepup + \stepup);
    }   
    
\end{tikzpicture}

$\sigma(H)$
\end{minipage}%
\hfill
\begin{minipage}[t]{.24\textwidth}
\centering
% \sigma^{2}(G)
\begin{tikzpicture}[rotate=270]
    \newcommand{\n}{4}                          %size of the permutations
    \newcommand{\width}{3}                       
    \newcommand{\height}{2}                       

    \pgfmathsetmacro\N{6}
    \pgfmathsetmacro\stepright{\width / \N}
    \pgfmathsetmacro\stepup{\height / \n}

%draw colors
\foreach\i/\a/\b/\c/\d in                       %this needs to be expanded if \n is bigger
{1/3/4/1/2,
2/2/4/1/3,
3/2/4/3/1,
4/1/4/3/2,
5/1/4/2/3,
6/3/4/2/1
}
{
    \newcommand{\sethgt}[1]
    {
        \if##1\a\def\hgt{1}\fi
        \if##1\b\def\hgt{2}\fi
        \if##1\c\def\hgt{3}\fi
        \if##1\d\def\hgt{4}\fi
    }
        \foreach\y in {1,...,\n}
        {
            \sethgt{\y} %sets \hgt
            \setclr{\y} %sets \clr
            \fill[\clr] (\i * \stepright, \hgt * \stepup) -- (\i * \stepright + \stepright, \hgt * \stepup) -- (\i * \stepright + \stepright, \hgt * \stepup + \stepup) -- (\i * \stepright, \hgt * \stepup + \stepup) -- cycle;
        }
}

%grid
    \foreach\y in {0,...,\n}
    {
        \draw (\stepright, \y * \stepup + \stepup) -- (\N * \stepright + \stepright, \y * \stepup + \stepup);
    }
    \foreach\x in {0,...,\N}
    {
        \draw (\x * \stepright + \stepright, \stepup) -- (\x * \stepright + \stepright, \n * \stepup + \stepup);
    }   
    
\end{tikzpicture}

$\sigma^{2}(G)$
\end{minipage}%
\hfill
\begin{minipage}[t]{.24\textwidth}
\centering
% \sigma^{3}(H)
\begin{tikzpicture}[rotate=270]
    \newcommand{\n}{4}                          %size of the permutations
    \newcommand{\width}{3}                       
    \newcommand{\height}{2}                       

    \pgfmathsetmacro\N{6}
    \pgfmathsetmacro\stepright{\width / \N}
    \pgfmathsetmacro\stepup{\height / \n}

%draw colors
\foreach\i/\a/\b/\c/\d in                       %this needs to be expanded if \n is bigger
{1/3/2/4/1,
2/2/3/4/1,
3/1/3/4/2,
4/3/1/4/2,
5/2/1/4/3,
6/1/2/4/3
}
{
    \newcommand{\sethgt}[1]
    {
        \if##1\a\def\hgt{1}\fi
        \if##1\b\def\hgt{2}\fi
        \if##1\c\def\hgt{3}\fi
        \if##1\d\def\hgt{4}\fi
    }
        \foreach\y in {1,...,\n}
        {
            \sethgt{\y} %sets \hgt
            \setclr{\y} %sets \clr
            \fill[\clr] (\i * \stepright, \hgt * \stepup) -- (\i * \stepright + \stepright, \hgt * \stepup) -- (\i * \stepright + \stepright, \hgt * \stepup + \stepup) -- (\i * \stepright, \hgt * \stepup + \stepup) -- cycle;
        }
}

%grid
    \foreach\y in {0,...,\n}
    {
        \draw (\stepright, \y * \stepup + \stepup) -- (\N * \stepright + \stepright, \y * \stepup + \stepup);
    }
    \foreach\x in {0,...,\N}
    {
        \draw (\x * \stepright + \stepright, \stepup) -- (\x * \stepright + \stepright, \n * \stepup + \stepup);
    }   
    
\end{tikzpicture}

$\sigma^{3}(H)$
\end{minipage}%

\caption{Visualization of the balanced Gray code for $\Pcadj_4$ as constructed in \Cref{thm:cyclicallyballanced} for $n=4$ where $1=$~purple, $2=$~blue, $3=$~green and $4=$~yellow and $\sigma := (1\,2\,3\,4)$.}
\label{fig:balancedAdjacentExample}
\end{figure}

\subsection{Small Rainbow Cycles and Large Rainbow Paths on the Permutahedron} \label{sec:permutahedron}

    \permutahedron*
    \begin{proof}
        Let us begin with $r=2$. Here are the flips of a 2-rainbow cycle in $P_5$:
        \[
            [(1,2),(2,3),(4,5),(2,3),(3,4),(1,2),(3,4),(4,5)].
        \]
        We now make the inductive step $n\to n+1$. By the induction hypothesis there is a 2-rainbow cycle $C$ in $P_n$. Obtain the path $P$ in $P_n$ by removing one of the occurrences of $(1,2)$ from $C$ and starting after the removed transposition. Note that we have defined the cycle as a list of transpositions, so that we obtain the actual cycle in $\Padj_n$ by choosing an initial permutation and then applying the transpositions as given. Now the 2-rainbow cycle in $P_{n+1}$ is obtained as follows:
        \begin{align*}
            &1,2,3,\ldots ,n-1,\underline{n,n+1}\\
            \to \qquad &\underline{1,\underline{2,3,\ldots,n-1,n+1}},n\\
            \underbrace{\to \ldots \to}_{P} \qquad &2,1,3,\ldots,n-1,\underline{n+1,n}\\
            \to \qquad &\underline{2,1},3,\ldots,n-1,n,n+1\\
            \to \qquad &1,2,3,\ldots,n-1,n,n+1.
        \end{align*}
        Here the underlined tuples are transposed to obtain the next permutation. Note that $P$ uses each adjacent transposition in $\{(i,i+1)\mid i\in[n]\}$ exactly twice except for $(1,2)$.\\
        Now we show that $P_n$ has a 3-rainbow cycle for $n\ge 3$ and $n$ odd. Here are the flips of a 3-rainbow cycle in $P_3$:
        \[
            [(1,2),(2,3),(1,2),(2,3),(1,2),(2,3)].
        \]
        We now make the inductive step $n\to n+2$. By the induction hypothesis there is a 3-rainbow cycle $C$ in $P_n$. Obtain the path $P$ in $P_n$ by removing one of the occurences of $(1,2)$ from $C$ and starting after the removed transposition. Now the 3-rainbow cycle in $P_{n+2}$ is obtained as follows:
        \begin{align*}
            &1,2,3,\ldots,n-1,n,\underline{n+1,n+2}\\
            \to\qquad &1,2,3,\ldots,n-1,\underline{n,n+2},n+1\\
            \to\qquad &\underline{\underline{1,\underline{2,3,\ldots,n-1,n+2}}},n,n+1\\
            \underbrace{\to \ldots \to}_{P}\qquad &2,1,3,\ldots,n-1,n+2,\underline{n,n+1}\\
            \to\qquad &\underline{2,1},3,\ldots,n-1,n+2,n+1,n\\
            \to\qquad &1,2,3,\ldots,n-1,\underline{n+2,n+1},n\\
            \to\qquad &1,2,3,\ldots,n-1,n+1,\underline{n+2,n}\\
            \to\qquad &1,2,3,\ldots,n-1,\underline{n+1,n},n+2\\
            \to\qquad &1,2,3,\ldots,n-1,n,n+1,n+2.
        \end{align*}
        In both constructions above it is obvious that all permutations obtained are distinct. Hence these really are cycles.
    \end{proof}

\section{Conclusion and Open Problems}\label{sec:conclusion}

\subsection{The Associahedron}

Although we have defined $r$-rainbow cycles in $\A_n$ with $r\sim 2^n$, there is still a lot of room at the top. Roughly speaking, the number of triangulations of the $n$-gon is on the order of $4^n$. A balanced Gray code for the Associahedron does not exist in most cases simply because the exact number of triangulations is not divisible by the number of possible diagonals. However, we believe that much larger rainbow cycles should exist.

A related open problem is the construction of Hamilton-cycles in $\A_n$ in which the flips are as balanced as possible.

\subsection{Permutation Gray Codes}

Regarding the generation of permutations, we have shown that $\Pall_n$ and $\Pcadj_n$ admit balanced Gray codes for all $n\in\mathbb{N}$, settling and advancing questions by Felsner, Gregor, Kleist, Merino, Mütze and Sering. In all our constructions major usage is made of cyclic rotations. As a result, our methods, at least out of the box, do not carry over to the permutahedron, here called $\Padj$, since they are not closed under adjacent transpositions. It is known that $\Padj_1$, $\Padj_3$ and $\Padj_5$ admit balanced Gray codes, however $\Padj_4$ does not \cite{Muetze2023}. Given our results, we restate the question (\cite[$P62$]{Muetze2023}) whether~$\Padj_n$ contains a balanced Gray code whenever $n$ is odd. Complementary, it would be interesting to find the maximal $r$-rainbow cycle in $\Padj_n$ for even $n$. Our constructions of $2$ and $3$-rainbow cycles for $\Padj_n$ present first ideas towards progress in this direction.

Furthermore, we believe that, from an algorithmic perspective, the creation of greedy algorithms for the generation of the introduced balanced codes for $\Pall_n$ and $\Pcadj_n$ is a worthwhile endeavor. It should be possible to choose the next permutation on the fly at each point during the generation without knowing any completed balanced Gray code for lower orders. This is because, in principle, any next permutation in the presented Gray codes does not depend on future positions in the used Gray codes for any lower order. If the next permutation can always be determined efficiently, then the Gray code can be efficiently generated.

\begin{figure}[h!]
    \centering
    \includegraphics[trim=10bp 10bp 10bp 10bp]{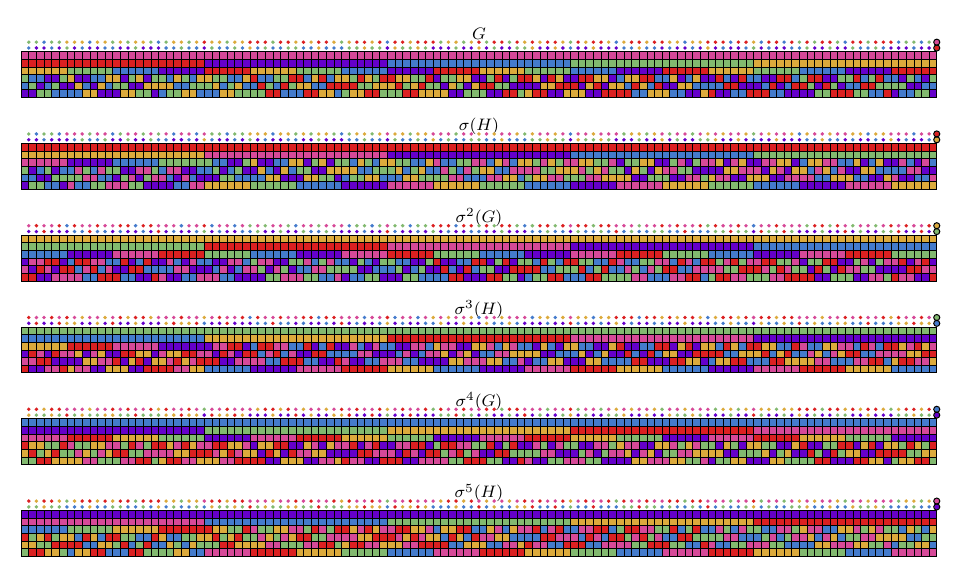}
    \caption{The balanced Gray code for $\Pall_6$ as constructed in the proof of \Cref{thm:perm_rainbow}.}
    \label{fig:Pall6}
\end{figure}

\begin{figure}[h!]
    \centering
    
\includegraphics[trim=10bp 10bp 10bp 10bp]{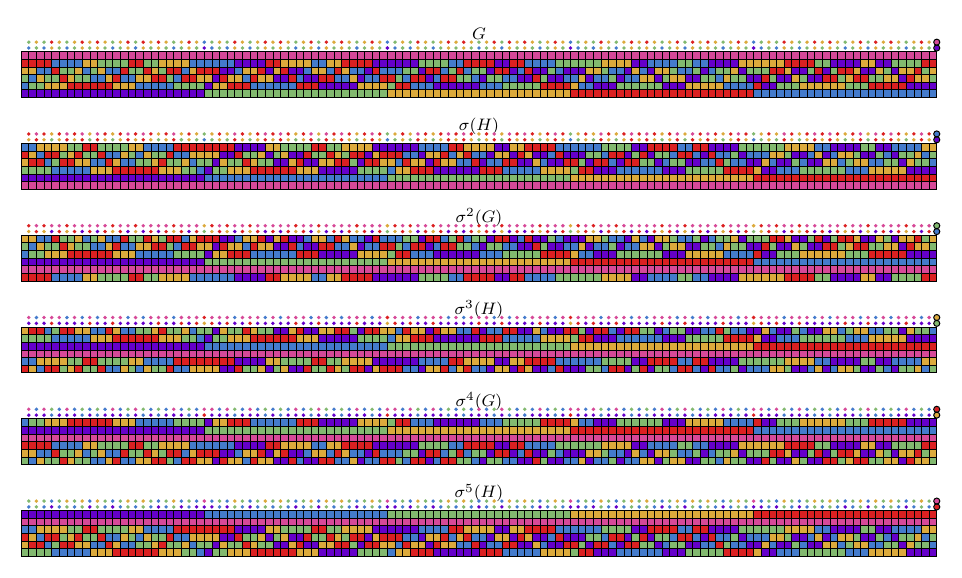}
    \caption{The balanced Gray code for $\Pcadj_6$ as constructed in the proof of \Cref{thm:cyclicallyballanced}.}
    \label{fig:Pcadj6}
\end{figure}

\subparagraph{Acknowledgements}
We wish to thank Torsten Mütze and Stefan Felsner for helpful discussions on this subject.

\printbibliography

\end{document}